\documentclass[journal,twoside]{IEEEtran}

\usepackage{amsmath,amssymb,amsthm,epsfig,dsfont,color,subfigure,empheq,enumerate}
\usepackage{stfloats,hyperref}
\usepackage{algorithmic,algorithm}

\DeclareMathOperator{\rank}{rank}

\DeclareMathOperator{\nullspace}{null}

\DeclareMathOperator{\diag}{diag}

\DeclareMathOperator*{\argmin}{argmin}

\theoremstyle{remark}

\input{mysymbol.sty}

\allowdisplaybreaks[4]

\begin{document}
\bstctlcite{IEEEexample:BSTcontrol}
\title{Monitoring and Optimization for Power Grids:\\
A Signal Processing Perspective$^{\dagger}$}
\author{Georgios B. Giannakis$^{\star}$,~\IEEEmembership{Fellow,~IEEE},
Vassilis Kekatos,~\IEEEmembership{Member,~IEEE}, Nikolaos
Gatsis,~\IEEEmembership{Member,~IEEE}, Seung-Jun Kim,~\IEEEmembership{Senior
Member,~IEEE}, Hao Zhu,~\IEEEmembership{Member,~IEEE},\\ and Bruce F.
Wollenberg,~\IEEEmembership{Fellow,~IEEE}
\thanks{$^{\dagger}$ Manuscript received August 31, 2012; revised December 20, 2012. Work in this paper was supported by the Inst. of Renewable Energy and the Environment (IREE) under grant no. RL-0010-13, Univ. of Minnesota.}
\thanks{$^{\star}$G. B. Giannakis, V. Kekatos, N. Gatsis, S.-J. Kim, and B. F. Wollenberg are with the Dept. of ECE and the Digital Technology Center (DTC) at the University of Minnesota, Minneapolis, MN 55455, USA. H. Zhu is with the Dept. of ECE at the University of Illinois, Urbana-Champaign. Emails:\{georgios,kekatos,gatsisn,seungjun,zhuh,wollenbe\}@umn.edu}}


\maketitle


\section{Introduction}\label{sec:intro}
Albeit the North American power grid has been recognized as the most
important engineering achievement of the 20th century, the modern power grid
faces major challenges~\cite{nae-report}. Increasingly complex
interconnections even at the continent size render prevention of the rare yet
catastrophic cascade failures a strenuous concern. Environmental incentives
require carefully revisiting how electrical power is generated, transmitted,
and consumed, with particular emphasis on the integration of renewable energy
resources. Pervasive use of digital technology in grid operation demands
resiliency against physical and cyber attacks on the power infrastructure.
Enhancing grid efficiency without compromising stability and quality in the
face of deregulation is imperative. Soliciting consumer participation and
exploring new business opportunities facilitated by the intelligent grid
infrastructure hold a great economic potential.

The smart grid vision aspires to address such challenges by capitalizing on
state-of-the-art information technologies in sensing, control, communication,
and machine learning \cite{AmW05,wollen11acc}. The resultant grid is
envisioned to have an unprecedented level of situational awareness and
controllability over its services and infrastructure to provide fast and
accurate diagnosis/prognosis, operation resiliency upon contingencies and
malicious attacks, as well as seamless integration of distributed energy
resources.

\subsection{Basic Elements of the Smart Grid}
A cornerstone of the smart grid is the advanced monitorability on its
assets and operations. Increasingly pervasive installation of the phasor
measurement units (PMUs) allows the so-termed synchrophasor
measurements to be taken roughly 100 times faster than the legacy supervisory
control and data acquisition (SCADA) measurements, time-stamped using the global positioning system (GPS) signals to capture the grid dynamics. In addition, the availability of
low-latency two-way communication networks will pave the way to
high-precision real-time grid state estimation and detection, remedial
actions upon network instability, and accurate risk analysis and post-event
assessment for failure prevention.

The provision of such enhanced monitoring and communication capabilities lays
the foundation for various grid control and optimization components. Demand
response (DR) aims to adapt the end-user power usage in response to energy
pricing, which is advantageously controlled by utility companies via smart meters~\cite{Hamilton-dr-pmag}. Renewable sources such as
solar, wind, and tidal, and electric vehicles are important pieces of the
future grid landscape. Microgrids will become widespread based on distributed
energy sources that include distributed generation and storage systems. Bidirectional power flow to/from the grid due to such
distributed sources has potentials to improve the grid economy and
robustness. New services and businesses will be generated through open grid
architectures and markets.

\subsection{SP for the Grid in a Nutshell: Past, Present, and Future}
Power engineers in the 60's were facing the problem of computing voltages at
critical points of the transmission grid, based on power flow readings taken
at current and voltage transformers. Local personnel manually collected these
readings and forwarded them by phone to a control center, where a set of
equations dictated by Kirchoff's and Ohm's laws were solved for the electric
circuit model of the grid. However, due to timing misalignment,
instrumentation inaccuracy, and modeling uncertainties present in these
measurements, the equations were always infeasible. Schweppe and others
offered a statistical signal processing (SP) problem formulation, and
advocated a least-squares approach for solving it~\cite{Schweppe70}---what
enabled the power grid monitoring infrastructure used pretty much invariant
till now \cite{Mo00}, \cite{AburExpositoBook}.

This is a simple but striking example of how SP expertise can have a strong
impact in power grid operation. Moving from the early 70's to nowadays, the
environment of the power system operation has become considerably more
complex. New opportunities have emerged in the smart grid context, necessitating a fresh look. As will be surveyed in this article, modern grid
challenges urge for innovative solutions that tap into diverse SP techniques
from estimation, machine learning, and network science.

Avenues where significant contribution can be made include power system state estimation (PSSE) in various renditions, as well as ``bad data'' detection and removal. As costly large-scale blackouts can be caused by rather minor outages in distant parts of the network, wide-area monitoring of the grid turns out to be a challenging yet essential goal~\cite{WAMS_PROC}. Opportunities abound in synchrophasor technology, ranging from judicious placement of PMUs to their role in enhancing observability, estimation accuracy, and bad data diagnosis. Unveiling topological changes given a limited set of power meter readings is a critical yet demanding task. Applications of machine learning to the power grid for clustering, topology inference, and Big Data processing 
for e.g., load/price forecasting constitute additional promising directions.

Power grid operations that can benefit from the SP
expertise include also traditional operations such as economic dispatch, power
flow, and unit commitment~\cite{WollenbergBook}, \cite{ShYaLi02},
\cite{ExpConCanBook}, as well as contemporary ones related to demand
scheduling, control of plug-in electric vehicles, and integration of
renewables. Consideration of distributed coordination of the partaking
entities along with the associated signaling practices and architectures
require careful studies by the SP, control, and optimization experts.

Without any doubt, computationally intelligent approaches based on SP
methodologies will play a crucial role in this exciting endeavor. From grid
informatics to inference for monitoring and optimization tools,
energy-related issues offer a fertile ground for SP growth whose time has
come.

The rest of the article is organized as follows. Modeling preliminaries for
power system analysis are provided in Sec.~\ref{sec:modeling}.
Sec.~\ref{sec:monitoring} deals with the monitoring aspect, delineating
various SP-intensive topics including state estimation and PMUs, as well as
the inference, learning and cyber-security tasks. Section~\ref{sec:control} is
devoted to grid optimization issues, touching upon both traditional problems
in economic power system operations, as well as more contemporary topics such
as demand response, electric vehicles, and renewables. The article is wrapped
up with a few open research directions in Sec.~\ref{sec:open}.

\section{Modeling Preliminaries}\label{sec:modeling}
Power systems can be thought of as electric circuits of even
continent-wide dimensions. They obey multivariate versions of
Kirchoff's and Ohm's laws, which in this section are overviewed
using a matrix-vector notation. As the focus is laid on alternating
current (AC) circuits, all electrical quantities involved (voltage,
current, impedance, power) are complex-valued. Further, quantities
are measured in the per unit (p.u.) system, which means that they
are assumed properly normalized. For example, if the ``base voltage"
is $138$ kV, then a bus voltage of $140$ kV is 1.01 p.u. The p.u.
system enables uniform single- and three-phase system analysis,
bounds the dynamic range of calculations, and allows for uniform
treatment over the different voltage levels present in the power
grid \cite{WollenbergBook}, \cite{ExpConCanBook}.

\begin{figure}
\centering
\includegraphics[width=0.95\linewidth]{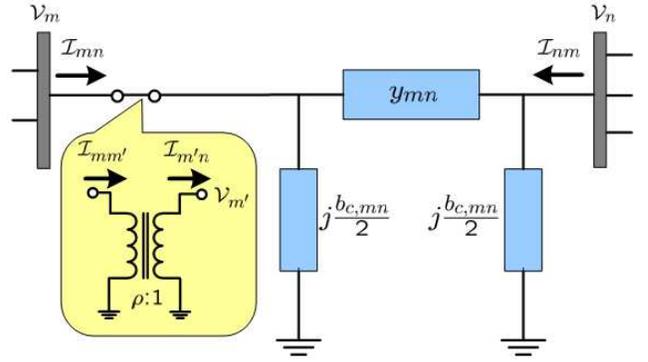}
\vspace*{-0.5em}
\caption{Equivalent $\pi$ model for a transmission line; yellow box
when an ideal transformer is also present [cf. \eqref{eq:i_mn2}].}
\label{fig:pi_model}
\vspace*{-1em}
\end{figure}

Consider first a power system module of two nodes, $m$ and $n$,
connected through a line. A node, also referred to as a bus in the
power engineering nomenclature, can represent, e.g., a generator or
a load substation. A line (a.k.a. branch) can stand for a
transmission or distribution line (overhead/underground), or even a
transformer. Two-node connections can be represented by the
equivalent $\pi$ model depicted in
Fig.~\ref{fig:pi_model}~\cite{MATPOWER,Mo00}, which entails the line
series impedance $z_{mn}:=1/y_{mn}$ and the total charging
susceptance $b_{c,mn}$. The former comprises a resistive part
$r_{mn}$ and a reactive (actually inductive) one $x_{mn}>0$, that is
$z_{mn}=r_{mn}+jx_{mn}$. The line series admittance
$y_{mn}:=1/z_{mn}=g_{mn}+jb_{mn}$ is often used in place of the
impedance. Its real and imaginary parts are called conductance and
susceptance, respectively. Letting $\mathcal{V}_m$ denote the
complex voltage at node $m$, $\mathcal{I}_{mn}$ the current flowing
from node $m$ to $n$, and invoking Ohm's and Kirchoff's laws on the
circuit of Fig.~\ref{fig:pi_model}, yields
\begin{equation}\label{eq:i_mn}
\mathcal{I}_{mn}=(jb_{c,mn}/2+y_{mn})\mathcal{V}_m-y_{mn}\mathcal{V}_n \;.
\end{equation}
The reverse-direction current $\mathcal{I}_{nm}$ is expressed symmetrically.
Unless $b_{c,mn}$ is zero, it holds that
$\mathcal{I}_{nm}\neq-\mathcal{I}_{mn}$. A small shunt susceptance $b_{s,mm}$
is typically assumed between every node $m$ and the ground (neutral),
yielding the current $\mathcal{I}_{mm}=jb_{s,mm}\mathcal{V}_m$.

Building on the two-node module, consider next a power system
consisting of a set $\mathcal{N}$ of $N_b$ buses along with a set
$\mathcal{E}$ of $N_l$ transmission lines. By Kirchoff's current
law, the complex current at bus $m$ denoted by $\mathcal{I}_m$ must
equal the sum of currents on the lines incident to bus $m$; that is,
\begin{align}
\mathcal{I}_{m} & =\sum_{n\in\mathcal{N}_m}\mathcal{I}_{mn}+\mathcal{I}_{mm} \notag\\
& =\left(\sum_{n\in\mathcal{N}_{m}}y_{mn} + y_{mm}\right) \mathcal{V}_{m}
- \sum_{n\in\mathcal{N}_{m}}y_{mn}\mathcal{V}_{n} \label{eq:i_m}
\end{align}
where $\mathcal{N}_m$ is the set of buses directly connected to bus
$m$, and
$y_{mm}:=j\left(b_{s,mm}+\sum_{n\in\mathcal{N}_m}b_{c,mn}/2\right)
:= j b_{mm}$. Collecting node voltages (currents) in the $N_b\times
1$ vector $\mathbf{v}$ ($\mathbf{i}$), leads to the multivariate
Ohm's law
\begin{equation}\label{eq:iYv}
\mathbf{i}=\mathbf{Y}\mathbf{v}
\end{equation}
where $\mathbf{Y}\in\mathbb{C}^{N_b\times N_b}$ is the so-termed
\emph{bus admittance matrix} with $(m,m)$-th diagonal entry
$\sum_{n\in\mathcal{N}_{m}}y_{mn} + y_{mm}$ and $(m,n)$-th
off-diagonal entry $-y_{mn}$ if $n\in\mathcal{N}_m$, and zero
otherwise (cf. \eqref{eq:i_m}). Matrix $\mathbf{Y}$ is symmetric and
more importantly sparse, thus facilitating efficient storage and
computations. On the contrary, the bus impedance matrix
$\mathbf{Z}$, defined as the inverse of $\mathbf{Y}$ (and not as the
matrix of bus pair impedances), is full and therefore it is seldom
used.

A major implication of \eqref{eq:iYv} is control of power flows. Let
$\mathcal{S}_m:=P_m+jQ_m$ be the complex power injected at bus $m$
whose real and imaginary parts are the active (reactive) power $P_m$
($Q_m$). Physically, $\mathcal{S}_m$ represents the power generated
and/or consumed by plants and loads residing at bus $m$. For bus $m$
and with $^*$ denoting conjugation, it holds that
$\mathcal{S}_m=\mathcal{V}_m\mathcal{I}_m^*$, or after collecting
all power injections in $\mathbf{s}\in\mathbb{C}^{N_b}$
($\diag(\mathbf{v})$ denotes a diagonal matrix holding $\mathbf{v}$
on its diagonal) one arrives at (cf.~\eqref{eq:iYv})
\begin{equation}\label{eq:svYv}
\mathbf{s}=\diag(\mathbf{v})\mathbf{i}^*=\diag(\mathbf{v})\mathbf{Y}^*\mathbf{v}^*.
\end{equation}
Complex power flowing from bus $m$ to a neighboring bus $n$ is
similarly given by
\begin{equation}\label{eq:smn}
\mathcal{S}_{mn}=\mathcal{V}_m\mathcal{I}_{mn}^*.
\end{equation}
The ensuing analysis pertains mainly to nodal quantities. However,
line quantities such as line currents and power flows over lines can
be modeled accordingly using \eqref{eq:i_mn} and \eqref{eq:smn}.

Typically, the complex bus admittance matrix is written in
rectangular coordinates as $\mathbf{Y}=\mathbf{G}+j\mathbf{B}$. Two
options become available from~\eqref{eq:svYv}, depending on whether
the complex nodal voltages are expressed in polar or rectangular
forms. The polar representation $\mathcal{V}_m=V_m e^{j\theta_m}$
yields [cf. \eqref{eq:i_m}]
\begin{subequations}\label{eq:power_flow_polar}
\begin{align}
P_m&=\sum_{n=1}^{N_b} V_mV_n\left( G_{mn}\cos\theta_{mn} +
B_{mn}\sin\theta_{mn} \right)\label{eq:power_flow_polar_active}\\
Q_m&=\sum_{n=1}^{N_b} V_mV_n\left( G_{mn}\sin\theta_{mn} -
B_{mn}\cos\theta_{mn} \right)\label{eq:power_flow_polar_reactive}
\end{align}
\end{subequations}
where $\theta_{mn}:=\theta_m-\theta_n$ $\forall m$. Since $P_m$ and
$Q_m$ depend on phase differences $\{\theta_{mn}\}$, power
injections $\{\mathcal{S}_m\}$ are invariant to phase shifts of bus
voltages. This explains why a selected bus called the reference,
slack, or swing bus is conventionally assumed to have zero voltage
phase without loss of generality.

If $\mathbf{Y}$ is known, the $2N_b$ equations in \eqref{eq:power_flow_polar} involve the variables $\{P_m,Q_m,V_m,\theta_m\}_{m=1 }^{N_b}$. Among the $4N_b$ nodal variables, (i) the reference bus has fixed $(V_m,\theta_m)$; (ii) pairs $(P_m,V_m)$ are controlled at generator buses (and are thus termed PV buses); while, (iii) power demands $(P_m,Q_m)$ are predicted for load buses (also called PQ buses). Fixing these $2N_b$ variables and solving the non-linear equations~\eqref{eq:power_flow_polar} for the remaining ones constitutes the standard power flow problem \cite[Ch.~4]{WollenbergBook}. Algorithms for controlling PV buses and predicting load at PQ buses are presented in Sec.~\ref{subsec:opf} and Sec.~\ref{subsubsec:load_prediction}, respectively.

Pairs $(P_{mn},V_{mn})$ satisfying (approximately) power flow
equations paralleling~\eqref{eq:power_flow_polar} can be found in
\cite[Ch.~3]{ExpConCanBook}. Among the approximations of the latter
as well as \eqref{eq:power_flow_polar}, the so called \emph{DC
model} is reviewed next due to its importance in grid monitoring and
optimization. The DC model hinges on three assumptions:\\
\textbf{(A1)} The power network is purely inductive, which means
that $r_{mn}$ is negligible. In high-voltage transmission lines, the
ratio $x_{mn}/r_{mn}=-b_{mn}/g_{mn}$ is large enough so that
resistances can be ignored and the conductance part $\mathbf{G}$ of
$\mathbf{Y}$ can be approximated by zero; \\ \textbf{(A2)} In
regular power system conditions, the voltage phase differences
across directly connected buses are small; thus, $\theta_{mn}\simeq
0$ for every pair of neighboring buses $(m,n)$, and the
trigonometric functions in~\eqref{eq:power_flow_polar} are
approximated as $\sin\theta_{mn}\simeq \theta_m - \theta_n$ and
$\cos\theta_{mn}\simeq 1$; and \\ \textbf{(A3)} Due to typical
operating conditions, the magnitude of nodal voltages is
approximated by one p.u.

Under (A1)-(A3) and upon exploiting the structure of $\mathbf{B}$ (cf.
\eqref{eq:iYv}), the model in \eqref{eq:power_flow_polar} boils down to
\begin{subequations}\label{eq:power_flow_dc}
\begin{align}
P_m&= -\sum_{n\neq m}b_{mn} \left(\theta_m-\theta_n\right)\label{eq:power_flow_dc_active}\\
Q_m&= -b_{mm} - \sum_{n\neq m}b_{mn}\left(V_m-V_n\right)\label{eq:power_flow_dc_reactive}
\end{align}
\end{subequations}
where $b_{mn}= -1/x_{mn}$ is the susceptance of the $(m,n)$ branch,
and in deriving \eqref{eq:power_flow_dc}, approximation of nodal
voltage magnitudes to unity implies $V_mV_n\simeq 1$, yet
$V_m\left(V_m-V_n\right)\simeq V_m-V_n$.

The DC model \eqref{eq:power_flow_dc} entails {\it linear} equations
that are neatly decoupled: active powers depend only on voltage
phases, whereas reactive powers are solely expressible via voltage
magnitudes. Furthermore, the linear dependence is on voltage
differences. In fact, since $P_{mn}=-b_{mn}(\theta_m-\theta_n)$ and
$b_{mn}<0$, active power flows across lines from the larger- to the
smaller-voltage phase buses.

Consider now the active subproblem described by
\eqref{eq:power_flow_dc_active}. Stacking the nodal real power
injections in $\mathbf{p}\in\mathbb{R}^{N_b}$ and the nodal voltage
phases in $\boldsymbol{\theta}\in\mathbb{R}^{N_b}$, leads to
\begin{equation}\label{eq:pBtheta}
\mathbf{p}=\mathbf{B}_x\boldsymbol{\theta}
\end{equation}
where the symmetric $\mathbf{B}_x$ is defined similar to
$\mathbf{Y}$ by only accounting for reactances. Specifically,
$[\mathbf{B}_x]_{mm}:=\sum_{n\in\mathcal{N}_m}x_{mn}^{-1}$ for all
$m$, and $[\mathbf{B}_x]_{mn}:=-x_{mn}^{-1}$, if $(m,n)$ line
exists, and zero otherwise.

An alternative representation of $\mathbf{B}_x$ is presented next.
Define matrix $\mathbf{D}:=\diag\left(\{x_l^{-1}\}_{l \in
\mathcal{E}}\right)$, and the branch-bus $N_l\times N_b$ incidence
matrix $\mathbf{A}$, such that if its $l$-th row $\mathbf{a}_l^T$
corresponds to the $(m,n)$ branch, then $[\mathbf{a}_l]_m:=+1$,
$[\mathbf{a}_l]_n:=-1$, and zero elsewhere. Based on these
definitions, $\mathbf{B}_x = \mathbf{A}^T\mathbf{D}\mathbf{A}$ can
be viewed as a weighted Laplacian of the graph
$(\mathcal{N},\mathcal{E})$ describing the power network. This in
turn implies that $\mathbf{B}_x$ is positive semidefinite, and the
all-ones vector $\boldsymbol{1}$ lies in its null space. Further,
its rank is $(N_b-1)$ if and only if the power network is connected.
Since $\mathbf{B}_x\mathbf{1}=\mathbf{0}$, it follows that
$\mathbf{p}^T\mathbf{1}=0$; stated differently, the total active
power generated equals the active power consumed by all loads, since
resistive elements and incurred thermal losses are ignored.

As a trivia, the terminology {\em DC model} stems from the fact that
\eqref{eq:pBtheta} models the AC power system as a purely resistive
DC circuit by identifying the active powers, reactances, and the
voltage phases of the former to the currents, the resistances, and
the voltages of the latter.

Coming back to the exact power flow model of \eqref{eq:svYv},
consider now expressing nodal voltages in rectangular coordinates.
If $\mathcal{V}_m=V_{r,m}+jV_{i,m}$ for all buses, it follows that
\begin{subequations}\label{eq:power_flow_rect}
\begin{align}
P_m&=V_{r,m} \sum_{n=1}^{N_b} \left(V_{r,n}G_{mn}-V_{i,n}B_{mn}\right) \notag \\
& \phantom{=}\: +\: V_{i,m} \sum_{n=1}^{N_b} \left(V_{i,n}G_{mn}+V_{r,n}B_{mn}\right)
\label{eq:power_flow_rect_active}\\
Q_m &= V_{i,m} \sum_{n=1}^{N_b} \left(V_{r,n}G_{mn}-V_{i,n}B_{mn}\right) \notag\\
& \phantom{=} \: -\: V_{r,m} \sum_{n=1}^{N_b} \left(V_{i,n}G_{mn}+V_{r,n}B_{mn}\right).
\label{eq:power_flow_rect_reactive}
\end{align}
\end{subequations}
Based on \eqref{eq:power_flow_rect_active}
and~\eqref{eq:power_flow_rect_reactive}, it is clear that (re)active
power flows depend quadratically on the rectangular coordinates of
nodal voltages. Because \eqref{eq:power_flow_rect} is not amenable
to approximations invoked in deriving \eqref{eq:power_flow_polar},
the polar representation has been traditionally preferred over the
rectangular one.

Before closing this section, a few words are due on modeling
transformers that were not explicitly accounted so far. Upon adding
the circuit surrounded by the yellow square to the model of
Fig.~\ref{fig:pi_model}, the possibility of having a transformer on
a branch is considered in its most general setting
\cite{ExpConCanBook}, \cite{MATPOWER}. An ideal transformer residing
on the $(m,n)$ line at the $m$-th bus side yields
$\mathcal{V}_m=\mathcal{V}_{m'} \rho_{mn}$ and
$\mathcal{I}_{m'n}=\rho_{mn}^*\mathcal{I}_{mm'}$, where
$\rho_{mn}:=\tau_{mn} e^{j\alpha_{mn}}$ is its turn ratio. Hence,
\eqref{eq:i_mn} readily generalizes to
\begin{equation}\label{eq:i_mn2}
\left[\begin{array}{c}
\mathcal{I}_{mn}\\ \mathcal{I}_{nm}
\end{array}\right]=
\left[\begin{array}{cc}
\frac{y_{mn}+jb_{c,mn}/2}{|\rho_{mn}|^2} & -\frac{y_{mn}}{\rho_{mn}^*}\\
-\frac{y_{mn}}{\rho_{mn}} & y_{mn}+jb_{c,mn}/2
\end{array}\right]
\left[\begin{array}{c}
\mathcal{V}_{m}\\ \mathcal{V}_{n}
\end{array}\right].
\end{equation}
Using \eqref{eq:i_mn2} in lieu of \eqref{eq:i_mn}, a similar
analysis can be followed with the exception that in the presence of
phase shifters, the corresponding bus admittance matrix $\mathbf{Y}$
will not be symmetric. Note though that the DC model of
\eqref{eq:pBtheta} holds as is, since it ignores the effects of
transformers anyway.

The multivariate current-voltage law (cf.~\eqref{eq:iYv}), the power
flow equations (cf.~\eqref{eq:power_flow_polar} or
\eqref{eq:power_flow_rect}), along with their linear approximation
(cf.~\eqref{eq:pBtheta}) and generalization (cf.~\eqref{eq:i_mn2}),
will play instrumental roles in the grid monitoring, control, and
optimization tasks outlined in the ensuing sections.

\section{Grid Monitoring}\label{sec:monitoring}

In this section, SP tools and their roles in various grid monitoring
tasks are highlighted, encompassing state estimation with associated
observability and cyber-attack issues, synchrophasor measurements, as well as
intriguing inference and learning topics.

\subsection{Power System State Estimation}\label{subsec:SE}
Simple inspection of the equations in Section~\ref{sec:modeling} confirms that all nodal and line quantities become available if one knows the grid parameters $\{y_{mn}\}$, and all nodal voltages ${\cal V}_{mn}$ that constitute the system state. Power system state estimation (PSSE) is an important module in the supervisory control and data acquisition (SCADA) system for power grid operation. Apart from situational awareness, PSSE is essential in additional tasks, namely load forecasting, reliability analysis, the grid economic operations detailed in Sec.~\ref{sec:control}, network planning, and billing \cite[Ch.~4]{ExpConCanBook}. Building on Sec.~\ref{sec:modeling}, this section reviews conventional solutions and recent advances, as well as pertinent smart grid challenges and opportunities for PSSE.

\subsubsection{Static State Estimation}\label{subsubsec:static}
Meters installed across the grid continuously measure electric
quantities, and forward them every few seconds via remote terminal
units (RTUs) to the control center for grid monitoring. Due to imprecise time signaling and the SCADA scanning process, conventional metering cannot utilize phase information of the AC waveforms. Hence, legacy measurements involve (active/reactive) power injections and flows, as well as voltage and current magnitudes on specific grid points.
Given the SCADA measurements and assuming stationarity over a scanning cycle, the PSSE
module estimates the state, namely all complex nodal voltages collected in $\mathbf{v}$. Recall that according to the power flow models presented in Sec.~\ref{sec:modeling}, all grid quantities can be expressed in terms of $\mathbf{v}$.
Thus, the $M\times 1$ vector of SCADA measurements can be modeled as $\mathbf{z} = \mathbf{h}(\mathbf{v}) + \boldsymbol{\epsilon}$, where $\mathbf{h}(\cdot)$ is a properly defined vector-valued function, and $\boldsymbol{\epsilon}$ captures
measurement noise and modeling uncertainties. Upon prewhitening, $\boldsymbol{\epsilon}$ can be assumed standard Gaussian. The maximum-likelihood estimate (MLE) of $\mathbf{v}$ can be then simply expressed as the nonlinear least-squares (LS) estimate \begin{align}\label{eq:SE_wls}
\hhatbbv :=& \arg \min_\bbv \|\mathbf{z} - \mathbf{h}(\mathbf{v})\|_2^2.
\end{align}
Prior information, such as zero-injection buses ($P_m=Q_m=0$) and
feasible ranges (of $V_m$ and $\theta_m$), can be included as
constraints in \eqref{eq:SE_wls}. In any case, the optimization
problem is nonconvex. For example, when states are expressed in
rectangular coordinates, the functions in $\mathbf{h}(\cdot)$ are
quadratic; cf. \eqref{eq:power_flow_rect}. In general, PSSE falls
under the class of nonlinear LS problems, for which Gauss-Newton
iterations are known to offer the ``workhorse" solution
\cite[Ch.~2]{AburExpositoBook}. Specifically, upon expressing
$\mathbf{v}$ in polar coordinates, the quadratic
$\mathbf{h}(\mathbf{v})$ can be linearized using Taylor's expansion
around a starting point. The Gauss-Newton method hence approximates
the cost in \eqref{eq:SE_wls} with a linear LS one, and relies on
its minimizer to initialize the subsequent iteration. This iterative
procedure is closely related to gradient descent algorithms for
solving nonconvex problems, which are known to encounter two issues:
i) sensitivity to the initial guess; and ii) convergence concerns.
Without guaranteed convergence to the global optimum, existing
variants improve numerical stability of the matrix inversions per
iteration \cite{AburExpositoBook}. In a nutshell, the grand
challenge so far remains to develop a solver attaining or
approximating the \textit{global optimum} at
\textit{polynomial-time}.

Recently, a \emph{semidefinite relaxation} (SDR) approach has been
recognized to develop polynomial-time PSSE algorithms with the
potential to find a globally optimal solution~\cite{ZhuGia12}, \cite{ZhuGia-SGComm12}.
Challenged by the nonconvexity of \eqref{eq:SE_wls}, the measurement
model is reformulated as a linear function of the outer-product
matrix $\bbV:=\bbv\bbv^H$, where the state is now expressed in
rectangular coordinates. This allows reformulating~\eqref{eq:SE_wls}
to a semidefinite program (SDP) with the additional constraint
$\rank(\bbV)=1$. Dropping the nonconvex rank constraint to acquire a
convex SDP has been well-appreciated in signal processing and
communications; see e.g., \cite{LuMa10}. The SDR-based PSSE has been
shown to approximate well the global optimum, while it is possible
to further improve computational efficiency by exploiting the SDP
problem structure~\cite{ZhuGia12}.

\subsubsection{Dynamic State Estimation}\label{subsusec:dynamic}
As power systems evolve in time, dynamic PSSE is well motivated
thanks to its predictive ability emerging when additional temporal
information is available. In practice, it is challenged by both the
unknown dynamics and the requirement of real-time implementation.
While the latter could become tractable with (extended) Kalman
filtering (KF) techniques, it is more difficult to develop simple
state-space models to capture the power system dynamics.

There have been various proposals for state transition models in
order to perform the prediction step, mostly relying on a
quasi-steady state behavior; see \cite{RoVaDy90} for a review of the
main developments. One simplified and widely used model poses a
``random-walk'' behavior expressing the state in polar coordinates
per time slot $t$ as $\mathbf{v}(t+1) = \mathbf{v}(t) + \bbw(t)$,
where $\bbw(t)$ is zero-mean white Gaussian with a diagonal
covariance matrix estimated online~\cite{Mo00}.  A more
sophisticated dynamical model reads $\mathbf{v}(t+1) = \bbF(t) \mathbf{v}(t) + \bbe(t) + \bbw(t)$, where $\bbF(t)$ is a diagonal transition matrix and $\bbe(t)$
captures the process mismatch. Recently, a quasi-static
state model has been introduced to determine $\bbe(t)$ by
approximating first-order effects of load data~\cite{BlKrIl08}.

For the correction step, the extended KF (EKF) is commonly used via
linearizing the measurement model around the state predictor
\cite{Mo00,RoVaDy90}. To overcome the reduced accuracy of EKF
linearization, unscented KF (UKF) of higher complexity has been
reported in \cite{VaTe11}. Particle filtering may also be of
interest if its computational efficiency can be tolerated by the
real-time requirements of power systems.

\subsubsection{Distributed State Estimation}\label{subsubsec:distributed}
Parallel and distributed solvers were investigated early on
\cite{Schweppe70}. The motivation was primarily computational, even
though additional merits of coordination across adjacent control
areas were also recognized. In vertically integrated electricity
markets, each local utility estimated its own state and modeled the
rest of the system at boundary points using only local measurements.
Adjacent power systems were connected via \emph{tie lines}, which
were basically used in emergency situations, and PSSE was performed
locally with limited interaction among control centers.

Currently, the deregulation of energy markets has led to
continent-wide interconnections that are split into subnetworks
monitored by independent system operators (ISOs). Increasing
amount of power is transferred over multiple control areas, and tie
lines must be accurately monitored for reliability and accounting
\cite{ExpositoPROC11}. The ongoing penetration of renewables further
intensifies long-distance power transfers, while their intermittent
nature calls for frequent monitoring. Interconnection-level PSSE is
therefore a key factor for modernizing power grids. Even though
advanced instrumentation can provide precise and timely measurements
(cf. Sec.~\ref{subsec:PMU}), an interconnection could consist of
thousands of buses. The latter together with privacy policies deem
decentralized PSSE a pertinent solution.

\begin{figure}
\centering
\includegraphics[width=0.5\linewidth]{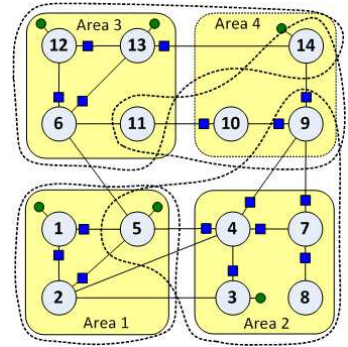}
\vspace*{-1em}
\caption{The IEEE 14-bus power system partitioned into four areas
\cite{PSTCA}. Dotted lassos show the buses belonging
to extended area states. PMU bus voltage (line current) measurements
are depicted by green circles (blue squares).}
\label{fig:ieee14partitioned}
\end{figure}

To understand the specifications of distributed PSSE, consider the
toy example of Fig.~\ref{fig:ieee14partitioned}. Area 2 consists of
buses $\{3,4,7,8\}$, but it also collects current measurements on
tie lines $\{(4,5),(4,9),(7,9)\}$. Its control center has two
options regarding these measurements: either to ignore them and
focus on the internal state, or to consider them and augment its
state by the external buses $\{5,9\}$. The first option is
statistically suboptimal; let alone it may incur observability loss
(check for example Area 3). For the second option, neighboring areas
should consent on shared variables. This way, agreement is achieved
over tie line charges and the global PSSE problem is optimally
solved.

It was early realized that for a chain of serially interconnected areas,
KF-type updates can be implemented incrementally in space \cite[Pt.
III]{Schweppe70}. For arbitrarily connected areas though, a two-level
approach with a global coordinator is required \cite{Schweppe70}: Local
measurements involving only local states are processed to estimate the
latter. Local estimates of shared states, their associated covariance
matrices, and tie line measurements are forwarded to a global coordinator.
The coordinator then updates the shared states and their statistics. Several
recent renditions of this hierarchical approach are available under the
assumption of local observability~\cite{ExpositoPROC11,MASEsurvey}. A central
coordinator becomes a single point of failure, while the sought algorithms
may be infeasible due to computational, communication, or policy limitations.
Decentralized solutions include block Jacobi iterations \cite{Conejo07}, and
the auxiliary problem principle~\cite{Ebrah00}. Local observability is waived
in \cite{XieChoiKar11}, where a copy of the entire high-dimensional state
vector is maintained per area, and linear convergence of the proposed
first-order algorithm scales unfavorably with the interconnection size. A
systematic framework based on the alternating direction method of multipliers
is put forth in~\cite{KeGi12}. Depending solely on existing PSSE software, it
respects privacy policies, exhibits low communication load, and its
convergence is guaranteed even in the absence of local observability.
Finally, for a survey on multi-area PSSE, refer to \cite{MASEsurvey}.

\begin{figure}
\centering
\subfigure[Bus/branch model.]{
\includegraphics[width=0.90\linewidth]{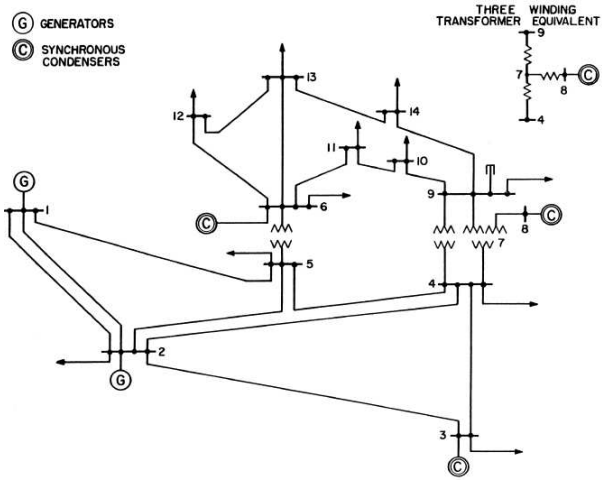}
\label{fig:ieee14}}\\
\subfigure[Bus section/switch model]{
\includegraphics[width=0.90\linewidth]{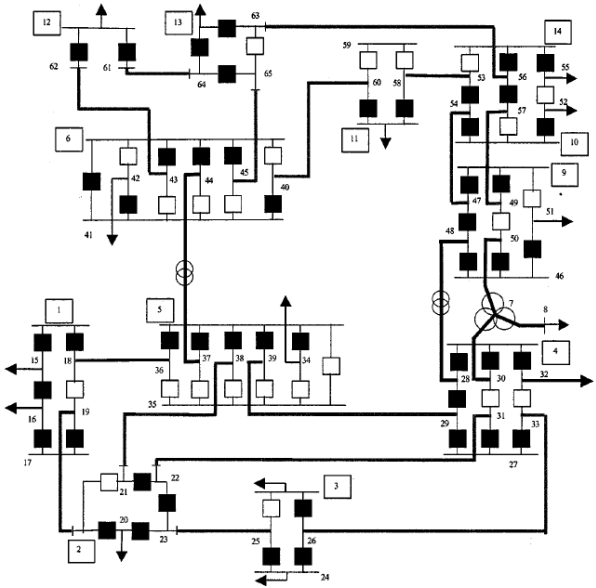}
\label{fig:ieee14expanded}} 
\caption{The IEEE 14-bus power system
benchmark \cite{PSTCA}: (a) The conventional model. (b) An assumed
substation-level model \cite{ExpVJ01}. Solid (hollow) squares
indicate closed (open) circuit breakers. The original 14 buses
preserve their numbering. Thick (thin) lines correspond to finite-
(zero-)impedance transmission lines (circuit breaker connections).}
\end{figure}

\subsubsection{Generalized State Estimation (G-SE)}\label{subsubsec:generalized}
PSSE presumes that grid connectivity and the electrical parameters involved
(e.g., line admittances) are known. Since these are oftentimes unavailable,
generalized state estimation (G-SE) extends the PSSE task to jointly
recovering them too \cite[Ch.~8]{AburExpositoBook},
\cite[Sec.~4.10]{ExpConCanBook}. PSSE operates on the bus/branch grid model;
cf.~Fig.~\ref{fig:ieee14}. A more meticulous view of this grid is offered by
the corresponding bus section/switch model depicted in
Fig.~\ref{fig:ieee14expanded}. This shows how a bus is partitioned by circuit breakers into sections (e.g., bus 1 to sections $\{1,15{-}19\}$), or how a
substation can appear as two different buses (e.g., sections $\{10,52{-}54\}$
and $\{14,55{-}57\}$ mapped to buses $10$ and $14$, respectively). Circuit
breakers are zero-impedance switching components and are used for seasonal,
maintenance, or emergency reconfiguration of substations. For some of them,
the status and/or the power they carry may be reported to the control center.
A topology processing unit collects this information and validates network
connectivity prior to PSSE \cite{Mo00}.

Even though topology malfunctions can be detected by large PSSE residual
errors, they are not easily identifiable \cite{AburExpositoBook}. Hence, joint PSSE with topology processing under the G-SE task has been a well-appreciated solution
\cite{Mo00}. G-SE essentially performs state estimation using the bus
section/switch model. Due to the zero impedances though, breaker flows are
appended to the system state. For regular transmission lines of unknown
status or parameters, G-SE augments the system state by their flows likewise.
In any case, to tackle the increased state dimensionality, breakers of known
status are treated as constraints: open (closed) breakers correspond to zero
flows (voltage drops). Practically, not all circuit breakers are monitored;
and even for those monitored, the reported status may be erroneous
\cite{AburExpositoBook}. Nowadays, G-SE is further challenged: the
penetration of renewables and DR programs will cause frequent substation 
reconfigurations. Yet, G-SE can be aided by advanced substation automation
and contemporary intelligent electronic devices (IEDs).

Identifying substation configuration errors has been traditionally treated by extending robust PSSE methods (cf.~Sec.~\ref{subsubsec:robust}) to the G-SE framework. Examples include the largest normalized residual test, and the least-absolute value and the Huber's estimators~\cite[Ch.~8]{AburExpositoBook}. To reduce the dimensionality of G-SE, an equivalent smaller-size model has been developed in \cite{ExpVJ01}. The method in \cite{KeGi12NAPS} leverages advances in compressive sampling and instrumentation technology. Upon regularizing the G-SE cost by $\ell_2$-norms of selected vectors, it promotes block sparsity on real and imaginary pairs of suspected breakers.

\subsection{Observability, Bad Data, and Cyber-attacks}\label{subsec:observability}
The PSSE module presumes that meters are sufficiently many and well distributed across the grid so that the power system is observable. Since this may not always be the case, observability analysis is the prerequisite of PSSE. Even when the set of measurements guarantees system state observability, resilience to erroneous readings should be solicited by robust PSSE methods. Nonetheless, specific readings (un)intentionally corrupted can harm PSSE results. This section studies these intertwined topics.

\subsubsection{Observability Analysis}\label{subsubsec:observability}
Given the network model and measurements, observability amounts to the ability of uniquely identifying the state $\mathbf{v}$. Even when the overall system is unobservable, power system operators are interested in observable islands. An observable island is a maximally connected sub-grid, whose states become observable upon selecting one of its buses as a reference. Identifying observable islands is important because it determines which line flows and nodal injections can be uniquely recovered. Identifying unobservable islands further provides candidate locations for additional (pseudo-)measurements needed to restore global observability. Pseudo-measurements are prior state information about e.g., scheduled generations, forecasted loads, or predicted values (based on historical data) to aid PSSE in the form of measurements with high-variance additive noise (estimation error).

Due to instrument failures, communication delays, and network reconfigurations, observability must be checked online. The analysis typically resorts to the DC model~\eqref{eq:power_flow_dc}, and hence, it can be performed separately per active and reactive subproblems thanks to the $P$-$\theta$ and $Q$-$V$ decoupling. Since power measurements oftentimes come in (re)active pairs, the observability results obtained for the active subproblem \eqref{eq:pBtheta} carry over to the reactive one, assuming additionally that at least one nodal voltage magnitude is available per observable island (the reactive analogue of the reference bus).

Commonly used observability checks include topological as well as numerical
ones; see \cite[Ch.~4]{AburExpositoBook} for a review. Topological
observability testing follows a graph-theoretic approach~\cite{ClKrDa83}.
Given the graph of the grid and the available set of measurements, this test
builds a maximal spanning tree. Its branches are either lines directly
metered or lines incident to a metered bus, while every branch should
correspond to a different measurement. If such a tree exists, the grid is
deemed observable; otherwise, the so-derived maximal spanning forest defines
the observable islands.

On the other hand, numerical observability considers the identifiability of the noiseless approximate DC model $\mathbf{z}=\mathbf{H}\boldsymbol{\theta}$~\cite{MoWu85A}. Linear system theory asserts that the state $\boldsymbol{\theta}$ is
observable if $\mathbf{H}$ is full column rank. Recall however that
active power measurements introduce a voltage phase shift ambiguity
(cf.~\eqref{eq:power_flow_polar}--\eqref{eq:power_flow_dc}). That is
why a power system with branch-bus incidence matrix $\mathbf{A}$ is
deemed observable simply if
$\mathbf{A}\boldsymbol{\theta}=\mathbf{0}$ for every
$\boldsymbol{\theta}$ satisfying
$\mathbf{H}\boldsymbol{\theta}=\mathbf{0}$, i.e.,
$\nullspace(\mathbf{H})\subseteq \nullspace(\mathbf{A})$. Observe now that the entries of $\mathbf{A}\boldsymbol{\theta}$ are proportional to line power flows. Hence, intuitively, whenever there is a non-zero power flow in the power
grid, at least one of its measurements should be non-zero for it to
be fully observable. When this condition does not hold, observable
islands can be identified via the iterative process developed in
\cite{MoWu85A}.

\subsubsection{Robust State Estimation by Cleansing Bad Data}\label{subsubsec:robust}
Observability analysis treats all measurements received as reliable and trustworthy. Nonetheless, time skews, communication failures, parameter uncertainty, and
infrequent instrument calibration can yield corrupted power system readings, also known as \emph{``bad data''} in the power engineering parlance. If bad data pass through simple screening tests, e.g., polarity or range checks, they can severely deteriorate PSSE performance. Coping with them draws methods from robust statistical SP to identify outlying measurements, or at least detect their presence in the measurement set.

Two statistical tests, namely the $\chi^2$-test and the largest normalized residual test (LNRT), were proposed in \cite[Part II]{Schweppe70}, and are traditionally used for bad data detection and identification, respectively \cite{Mo00}, \cite[Ch.~5]{AburExpositoBook}. Both tests rely on the model $\mathbf{z} = \mathbf{H}\boldsymbol{\theta}+\boldsymbol{\epsilon}$, assuming a full column rank $m\times n$ matrix $\mathbf{H}$ and a zero voltage phase at the reference bus. The two tests check the residual error of the LS estimator which can be expressed as $\mathbf{r}:=\mathbf{P}\mathbf{z}=\mathbf{P}\boldsymbol{\epsilon}$, where $\mathbf{P}:=\mathbf{I} - \mathbf{H}(\mathbf{H}^T\mathbf{H})^{-1}\mathbf{H}^T$ satisfying
$\mathbf{P}=\mathbf{P}^T=\mathbf{P}^2$. Apparently, when $\boldsymbol{\epsilon}$ is standardized Gaussian, $\mathbf{r}$ is Gaussian too with covariance $\mathbf{P}$; hence, $\|\mathbf{r}\|_2^2$ follows a $\chi^2$ distribution with $(m-n)$ degrees of freedom. The $\chi^2$-test then declares an LS-based PSSE possibly affected by outliers whenever $\|\mathbf{r}\|_2^2$ exceeds a predefined threshold.

LNRT exploits further the Gaussianity of $\mathbf{r}$. Indeed, as $r_i/\sqrt{P_{i,i}}$ should be standard Gaussian for all $i$ when bad data are absent, LNRT finds the maximum absolute value among these ratios and compares it against a threshold to identify a single bad datum \cite[Sec.~5.7]{AburExpositoBook}. Practically, if a bad datum is detected, it is removed from the measurement set, and the LS estimator is re-computed. The process is repeated till no bad data are identified. Successive LS estimates can be efficiently computed using recursive least-squares (RLS). The LNRT is essentially the leave-one-out approach, a classical technique for identifying single outliers. Interesting links between outlier identification and $\ell_0$-(pseudo)-norm minimization are presented in \cite{KoJi11} and \cite{KeGi12} under the Bayesian and the frequentist frameworks, respectively.

Apart from the two tests treating bad data a posteriori, outlier-robust estimators, such as the least-absolute deviation, the least median of squares, or Huber's estimator have been considered too; see \cite{AburExpositoBook}. Recently, $\ell_1$-norm based methods have been devised; see e.g., \cite{KoJi11}, \cite{XuWaTang11}, \cite{KeGi12}.

Unfortunately, all bad data cleansing techniques are vulnerable to the so called ``critical measurements'' \cite{AburExpositoBook}. A measurement is critical if once removed from the measurement set, the power system becomes unobservable. If for example one removes the current measurement on line $(7,8)$ from the grid of Fig.~\ref{fig:ieee14partitioned}, then bus $8$ voltage cannot be recovered. Actually, it can be shown that the $i$-th measurement is critical if the $i$-th column of $\mathbf{P}$ is zero, which translates to $r_i$ being always zero too. Due to the latter, the LNRT is undefined for critical measurements.


Intuitively, a critical measurement is the only observation related to some state. Thus, this measurement cannot be cross-validated or questioned as an outlier, but it should be blindly trusted. The existence of critical measurements in PSSE reveals the connection between bad data and observability analysis. Apparently, the notion of critical measurements can be generalized to multiple simultaneously corrupted readings. Even though such events are naturally rare, their study becomes timely nowadays under the threat of targeted cyber-attacks as explained next.

\subsubsection{Cyber-attacks}\label{subsubsec:cyber}
As a complex cyber-physical system spanning a large geographical area, the power grid inevitably faces challenges in terms of cyber-security. With more data acquisition and two-way communication required for the future grid, enhancing cyber-security is of paramount importance. From working experience in dealing with the Internet and telecommunication networks, there is potential for malicious and well-motivated adversaries to either physically attack the grid infrastructure, or remotely intrude the SCADA system. Among all targeted power grid monitoring and control operations, the PSSE task in Sec.~\ref{subsec:SE} appears to be of extreme interest as adversaries can readily mislead operators and manipulate electric markets by altering the system state \cite{KoJi11}, \cite{XiMo11}.

Most works analyzing cyber-attacks consider the linear measurement model modified as $\bbz =\bbH \bbtheta + \bbepsilon + \bba$, where the attack vector $\bba$ has non-zero entries corresponding to compromised meters. It was initially pointed out in \cite{LiNi11} that if the adversary knows $\mathbf{H}$, the attack $\bba$ can be constructed to lie in the range space of $\bbH$ so that the system operator can be arbitrarily misled. Under such a scenario, the attack cannot be detected. Such attacks are related to the observability and bad data analysis described earlier, since by deleting the rows of $\bbH$ corresponding to the nonzero entries of $\bba$, the resultant system becomes unobservable \cite{KoJi11}. Various strategies to construct $\bba$ have been derived in \cite{LiNi11}, constrained by the number of counterfeit meters; see also \cite{KoJi11} for the minimum number of such meters. Cyber-attacks under linear state-space models are considered in \cite{PasDorBul12}.

A major limitation of existing works lies in the linear measurement model assumption, not to mention the practicality of requiring attackers to know the full system configuration. Attacks in nonlinear measurement models for AC systems are studied in \cite{ZhGiaNAPS12}. Granted that a nonlinear PSSE model can be approximated around a given state point, it is not obvious how the attacker can acquire such dynamically varying information in real time in order to construct the approximation. This requires a per-adversary PSSE and assessment of a significant portion of meter measurements. On the defender's side, robustifying PSSE against bad data is a first countermeasure. Since cyber-attacks can be judiciously designed by adversaries, they may be more challenging to identify, thus requiring further prior information e.g., on the state vector statistics~\cite{KoJi11}.


\subsection{Phasor Measurement Units}\label{subsec:PMU}

\subsubsection{Phasor Estimation}\label{subsubsec:PMU_measurements}
PMUs are contemporary devices complementing legacy (SCADA) meters in
advancing power system applications via their high-accuracy and
time-synchronized measurements \cite{PhTh08}. Different from SCADA
meters which provide amplitude (power) related information, PMUs
offer also phase information. At the implementation level, current
and voltage transformers residing at substations provide the analog
input waveforms to a PMU. After anti-alias filtering, each one of
these analog signals is sampled at a rate several times the nominal
power system frequency $f_0$ (50/60~Hz). If the signal of interest
has frequency $f_0$, its phasor information (magnitude and phase)
can be obtained simply by correlating a window of its samples with
the sampled cosine and sine functions, or equivalently by keeping
the first (non-DC) discrete Fourier transform component. Such
correlations can be implemented also recursively. Since power system
components operate in the frequency range $f_0\pm 0.5$~Hz, acquiring
phasor information for off-nominal frequency signals has been also
considered \cite[Ch.~3]{PhTh08}.

The critical contribution of PMU technology to grid instrumentation
is time-tagging. Using precise GPS timing (the one pulse-per-second
signal), synchrophasors are time-stamped at the universal time
coordinated (UTC). PMU data can thus be consistently aggregated
across large geographic areas. Apart from phasors, PMUs acquire the
signal frequency and its frequency derivative too. Data from several
PMUs are collected by a phasor data concentrator (PDC) which
performs time-aligning, local cleansing of bad data, and potentially
data compression before forwarding data flows to the control center.
The IEEE standards C37.118.1/2-2011 determine PMU functional
requirements.

\subsubsection{PMU Placement}\label{subsubsec:PMU_placement}
Although PMU technology is sufficiently mature, PMU penetration has been
limited so far, mainly due to the installation and networking costs involved
\cite{WAMS_PROC}. Being the key technology towards wide area monitoring
though guarantees their wide deployment. During this instrumentation stage,
prioritizing PMU locations is currently an important issue for utilities and
reliability operators worldwide. Many PMU placement methods are based on the
notion of topological observability; cf. Sec.~\ref{subsubsec:observability}.
A search algorithm for placing a limited number of PMUs on a maximal spanning
forest is developed in \cite{NuPh05}.
Even though topological observability in general does not imply numerical
observability, for practical measurement matrices it does \cite{Mo00}. In any
case, a full column rank yet ill-conditioned linear regression matrix can
yield numerically unstable estimators. Estimation accuracy rather than
observability is probably a more meaningful criterion. Towards that end, PMU
placement is formulated as a variation of the optimal experimental design
problem in \cite{LiNeIl11}, \cite{KeGiWo12}. The approach in \cite{LiNeIl11}
considers estimating voltage phases only, ignores PMU current measurements,
and proposes a greedy algorithm. In \cite{KeGiWo12}, the state is expressed
in rectangular coordinates, all PMU measurements are considered, and the SDP
relaxation of the problem is solved via a projected gradient algorithm. For a
detailed review of PMU placements, the reader is referred to \cite{MaKoGe12}.

\subsubsection{State Estimation with PMUs}\label{subsubsec:PMU_SE}
As explained in Sec.~\ref{subsubsec:static}, PSSE is conventionally performed
using SCADA measurements \cite[Ch.~12]{WollenbergBook}. PMU-based PSSE
improves estimation accuracy when conventional and PMU measurements are
jointly used~\cite{PhThKa86,PhTh08}. However, aggregating conventional and
synchrophasor readings involves several issues. First, SCADA measurements are
available every 4 secs, whereas 30-60 synchrophasors can be reported per sec.
Second, explicitly including conventional measurements reduces the linear
PMU-based PSSE problem into a non-linear one. Third, compatibility to
existing PSSE software and phase alignment should be also considered. An
approach to address these challenges is treating SCADA-based estimates as
pseudo-measurements during PMU-driven state estimation \cite{PhTh08}.
Essentially, the slower rate SCADA-based state estimates, expressed in
rectangular coordinates, together with their associated covariance matrix can
be used as a Gaussian prior for the faster rate linear PSSE problem based on
PMU measurements \cite{PhTh08}, \cite{KeGiWo12}. Regarding phase alignment,
as already explained SCADA-based estimates assume the phase of the reference
bus to be zero, whereas PMUs record phases with respect to GPS timing.
Aligning the phases of the two estimates can be accomplished by
PMU-instrumenting the reference bus, and then simply adding its phase to all
SCADA-based state estimates \cite{PhTh08}.

Synchrophasor measurements do not contribute only to PSSE. Several
other monitoring, protection, and control tasks, ranging from local
to interconnection-wide scope can benefit from PMU technology.
Voltage stability, line parameter estimation, dynamic line rating,
oscillation and angular separation monitoring, small signal analysis
are just a few entries from the list of targeted applications
\cite{WAMS_PROC}, \cite{PhTh08}.

\subsection{Additional Inference and Learning Issues}\label{subsec:learning}
PSSE offers a prototype class of problems that SP tools can be readily
employed to advance grid monitoring performance, especially after leveraging
recent PMU technology to complement SCADA measurements. However, additional
areas can benefit from SP algorithms applied to change detection, estimation,
classification, prediction, and clustering aspects of the grid.

\subsubsection{Line Outage Identification}\label{subsec:topology}
Unexpected events, such as a breaker failure, a tree fall, or a lightning
strike, can make transmission lines inoperative. Unless the control center
becomes aware of the outage promptly, power generation and consumption will remain almost unchanged across the grid. Due to flow conservation though,
electric currents will be automatically altered in the outaged transmission network. Hence, shortly after, a few operating lines may exceed their ratings and
successively fail. A cascading failure can spread over interconnected systems
in a few minutes and eventually lead to a costly grid-wide blackout in less
than an hour. Timely identifying line outages, or more generally abrupt
changes in line parameters, is thus critical for wide-area monitoring.

\begin{figure}
\centering
\includegraphics[width=1.0\linewidth]{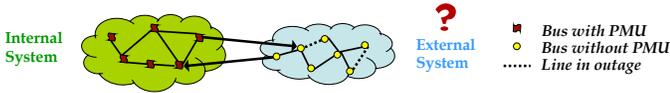}
\vspace*{-1em}
\caption{Internal system identifies line outages occurred in the
external system.} \label{fig:loid}
\vspace*{-1em}
\end{figure}

One could resort to the generalized PSSE module to identify line
outages (cf. Sec.~\ref{subsubsec:generalized}). Yet most existing
topology processors rely on data of the local control area (a.k.a.
\emph{internal}) system; see also Fig.~\ref{fig:loid}. On the other
hand, flow conservation can potentially reveal line changes even in
\emph{external} systems. This would be a non-issue if inter-system
data were available at a sufficiently high rate. Unfortunately, the
system data exchange (SDX) module of the North-American Electric
Reliability Corporation (NERC) can provide the grid-wide basecase
topology only on an hourly basis \cite{tate08tps}, while the
desideratum here is near-real-time line monitoring. In a nutshell,
each internal system needs to timely identify line changes even in
the external systems, relying only on local data and the
infrequently updated basecase topology.

To concretely lay out the problem, consider the pre- and post-event
states, and let $\tilde{\mathcal{E}} \subset \mathcal{E}$ denote the
subset of lines in outage. Suppose that the interconnected grid has
reached a stable post-event state, and it remains connected
\cite{tate08tps}. With reference to the linear DC model in
\eqref{eq:pBtheta}, its post-event counterpart reads $\bbp'=
\bbp+\bbeta = \bbB'_x \bbtheta'$, where $\bbeta$ captures small
zero-mean power injection perturbations. Recalling from
Sec.~\ref{sec:modeling} that $\mathbf{B}_x = \mathbf{A}^T \mathbf{D}
\mathbf{A}$, the difference $\tdbbB_x := \bbB_x - \bbB'_x$ can be
expressed as $\tilde{\mathbf{B}}_x=\sum_{\ell\in\tilde{\ccalE}}
x_\ell^{-1} \bba_\ell \bba_\ell^T$. With $\tdbbtheta:= \bbtheta' -
\bbtheta$, the ``difference model'' can be written as
$\mathbf{B}\tilde{\boldsymbol{\theta}} =
\sum_{\ell\in\tilde{\ccalE}} m_\ell \bba_\ell + \bbeta$, where
$m_\ell : = \bba_\ell^T \bbtheta'/x_\ell$, $\forall
\ell\in\tilde{\ccalE}$. Based on the latter, to identify
$\tilde{\ccalE}$ of a given cardinality $N_{l}^o:=|\tilde{\ccalE}|$,
one can enumerate all $\binom{N_l}{N_l^o}$ possible topologies in
outage, and select the one offering the minimum LS fit. Such an
approach incurs combinatorial complexity, and has thus limited the
existing exhaustive search methods to identifying single
\cite{tate08tps}, or at most double line outages \cite{tate09pesgm}.
A mixed-integer programming approach was proposed in
\cite{abur10naps}, which again deals with single line outages.

To bypass this combinatorial complexity, \cite{ZhGi12} considers an
overcomplete representation capturing all possible line outages. By
constructing an $N_l \times 1$ vector $\bbm$, whose $\ell$-th entry
equals $m_\ell$, if $\ell \in \tilde{\ccalE}$, and $0$ otherwise, it
is possible to reduce the previous model to a sparse linear
regression one given by
\begin{align}
\bbB_x \tdbbtheta= \bbA^T \bbm+ \bbeta. \label{smodel}
\end{align}
Since the control center only has estimates of the internal bus
phases, it is necessary to solve \eqref{smodel} for $\tdbbtheta$,
and extract the rows corresponding to the internal buses. This leads
to a linear model slightly different from \eqref{smodel}; but thanks
to the overcomplete representation, identifying $\tilde{\ccalE}$
amounts to recovering $\bbm$. The key point here is the small number
of line outages $(N_{l,o} \ll N_l)$ that makes the sought vector
$\bbm$ \emph{sparse}. Building on compressive sampling approaches,
sparse signal recovery algorithms have been tested in \cite{ZhGi12}
using IEEE benchmark systems, and near-optimal performance was
obtained at computational complexity growing only linearly in the
number of outages.

\subsubsection{Mode Estimation}\label{subsubsec:modal}
Oscillations emerge in power systems when generators are
interconnected for enhanced capacity and reliability. Generator
rotor oscillations are due to lack of damping torque, and give rise
to oscillations of bus voltages, frequency, and (re)active power
flows. Oscillations are characterized by the so-termed
electromechanical modes, whose properties include frequency,
damping, and shape \cite{Ku94}. Depending on the size of the power
system, modal frequencies are often in the range of $0.1-2$~Hz.
While a single generator usually leads to local oscillations at the
higher range ($1-2$~Hz), inter-area oscillations among groups of
generators lie in the lower range ($0.1-1$~Hz). Typically, the
latter ones are more troublesome, and without sufficient damping
they grow in magnitude and may finally result in even grid breakups.
Hence, estimating electromechanical modes, especially the
low-frequency ones, is truly important, and known as the
small-signal stability problem in power system analysis \cite{Ku94}.

\begin{figure}
\centering
\includegraphics[width=0.80\linewidth]{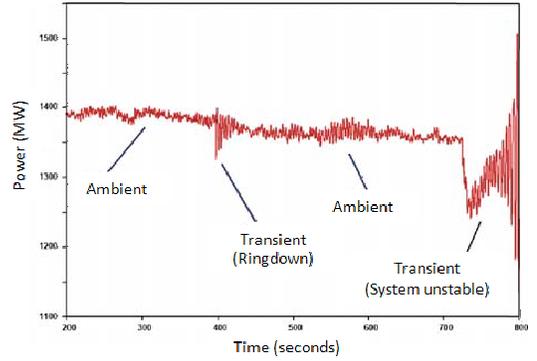}
\vspace*{-1em}
\caption{Real power flow on a major transmission line during the
1996 Western North American power system breakup~\cite{TrPi09}.}
\label{fig:modaldata}
\vspace*{-1em}
\end{figure}

Albeit near-and-dear to SP expertise on retrieving harmonics, modal
estimation is challenging primarily due to the nonlinear and time varying
properties of power systems, as well as the co-existence of several
oscillation modes at nearby frequencies. Fortunately, the system behaves
relatively linearly when operating at steady state, and can thus be
approximated by the continuous-time vector differential equations
${\dot{\bbx}}(t) = \bbA_x(t) \bbx(t) + \bbB_u \bbu(t) +\bbw(t)$, where the
eigenvalues of $\bbA_x(t)$ characterize the oscillation modes, and $\bbu(t)$
and $\bbw(t)$ correspond to the exogenous input and the random perturbing
noise, respectively. Assuming linear dynamic state models, mode estimation
approaches are either model- or measurement-based. The former construct the
exact nonlinear differential equations from system configurations, and then
linearize them at the steady-state to obtain $\bbA_x(t)$ for estimating
electromechanical modes~\cite{Mi05}. In measurement-based methods,
oscillation modes are acquired directly by peak-picking the spectral
estimates obtained using linear measurements $\bbx(t)$~\cite{TrPi09}. Since
the complexity of model-based methods grows with the network size,
scalability issues arise for larger systems. With PMUs, modes can be
estimated directly from synchrophasors, and even updated in real time.

Depending on the input $\bbu(t)$, the measurements are either ambient, or
ring-down (a.k.a. transient), or probing; see e.g., Fig.~\ref{fig:modaldata}.
With only random noise $\bbw(t)$ attributed to load perturbations, the system
operates under an equilibrium condition and the ambient measurements look
like pseudo-noise. A ring-down response occurs after some major disturbance,
such as line tripping or a pulse input $\bbu(t)$, and results in observable
oscillations. Probing measurements are obtained after intentionally injecting
known pseudo-random inputs (probing signals), and can be considered as a
special case of ring-down data. Missing entries and outliers are also
expected in meter measurements, hence robust schemes are of interest for mode
estimation~\cite{ZhTr08}. Measurement-based algorithms can be either batch or
recursive. In batch modal analysis, off-line ring-down data are modeled as a
sum of damped sinusoids and solved using e.g., Prony's method to obtain
linear transfer functions. Ambient data are handled by either parametric or
nonparametric spectral analysis methods \cite{TrPi09}. To recursively
incorporate incoming data, several adaptive SP methods have been successfully
applied, including least-mean squares (LMS) and RLS \cite{ZhTr08}. Apart from
utilizing powerful statistical SP tools for mode estimation, it is also
imperative to judiciously design efficient probing signals for improved
accuracy with minimal impact to power system operations \cite{TrPi09}.

\subsubsection{Load and Electricity Price Forecasting} \label{subsubsec:load_prediction}
Smooth operation of the grid depends heavily on load forecasts. Different
applications require load predictions of varying time scales. Minute- and hour-ahead load estimates are fed to the
unit commitment and economic dispatch modules as described in Sec.~\ref{subsec:opf}.
Predictions at the week scale are used for reliability purposes and
hydro-thermal coordination; while forecasts for years ahead facilitate
strategic generation and transmission planning. The granularity of load
forecasts varies spatially too, ranging from a substation, utility, to an
interconnection level. Load forecasting tools are essential for electricity
market participants and system operators. Even though such tools are widely
used in vertically organized utilities, balancing supply and demand at a
deregulated electricity market makes load forecasting even more important. At
the same time, the introduction of electric vehicles and DR programs further
complicates the problem.

Load prediction can be simply stated as the problem of inferring
future power demand given past observations. Oftentimes, historical
and predicted values of weather data (e.g., temperature and
humidity) are included as prediction variables too. The particular
characteristics of power consumption render it an intriguing
inference task. On top of a slowly increasing trend, load exhibits
hourly, weekly, and seasonal periodicities. Holidays, extreme
weather conditions, big events, or a factory interruption create
outlying data. Moreover, residential, commercial, and industrial
consumers exhibit different power profiles. Apart from the predicted
load, uncertainty descriptors such as confidence intervals are
important. Actually, for certain reliability and security
applications, daily, weekly, or seasonal peak values are critically
needed.

Several statistical inference methods have been applied for load
forecasting: ordinary linear regression; kernel-based regression and
support vector machines; time series analysis using auto-regressive
(integrated) moving average (with exogenous variables) models (ARMA,
ARIMA, ARIMAX); state-space models with Kalman and particle
filtering; neural networks, expert systems, and artificial
intelligence approaches. Recent academic works and current industry
practices are variations and combinations of these themes reviewed
in \cite[Ch.~2]{ShYaLi02}. Low-rank models for load imputation have been pursued in \cite{MaGia12}.

Load forecasting is not the only prediction task in modern power systems.
Under a deregulated power industry, market participants can also leverage
estimates of future electricity prices. To appreciate the value of such
estimates, consider a day-ahead market: an ISO determines the prices
of electric power scheduled for generation and consumption at the
transmission level during the 24 hours of the following day. The ISO collects
the hourly supply and demand bids submitted by generator owners and
utilities. Using the optimization methods described later in
Sec.~\ref{subsec:opf}, the grid is dispatched in the most economical way
while complying with network and reliability constraints. The output of this
dispatch are the power schedules for generators and utilities, along with
associated costs. Modern electricity markets are complex. Trading and hedging
strategies, weather and life patterns, fuel prices, government policies,
scheduled and random outages, reliability rules, all these factors influence
electricity prices. Even though prices are harder to predict than loads, the
task is truly critical in financial decision making \cite{AmHe06}. The
solutions proposed so far include econometric methods, physical system
modeling, time series and statistical methods, artificial intelligence
approaches, and kernel-based approaches; see e.g., \cite{AmHe06}, \cite{WuSh10}, \cite{KeLiVeGia13} and references therein.

\subsubsection{Grid Clustering}\label{subsubsec:clustering}

\begin{figure}
\centering
\includegraphics[width=1.0\linewidth]{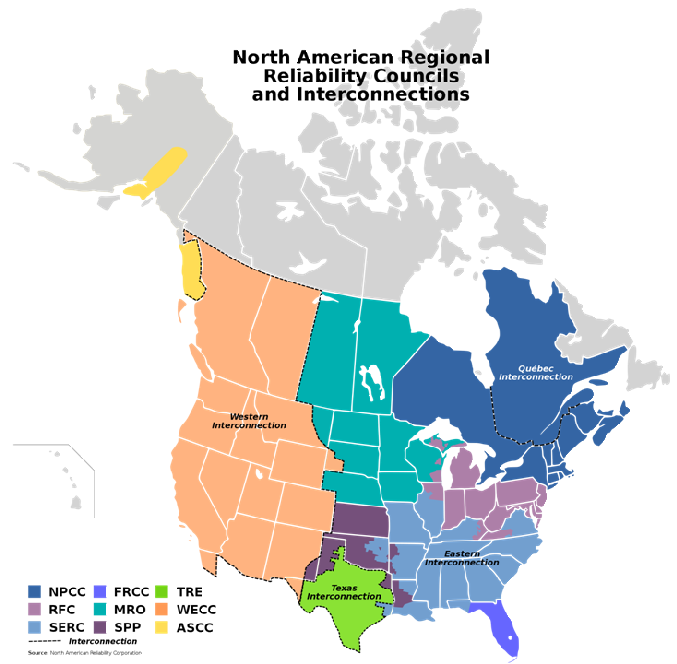}
\vspace*{-1em}
\caption{NERC's regional reliability councils and interconnections
[Source: \url{http://en.wikipedia.org/wiki/File:NERC-map-en.svg}]}
\label{fig:nerc}
\vspace*{-1em}
\end{figure}

Modularizing power networks is instrumental for grid operation as it
facilitates decentralized and parallel computation. Partitioning the
grid into control regions can also be beneficial for implementing
``self-healing" features, including islanding under contingencies
\cite{LiLi10}. For example, after catastrophic events, such as
earthquakes, alternative power supplies from different management
regions may be necessary due to power shortage and system
instability. Furthermore, grid partitioning is essential for the
zonal analysis of power systems, to aid load reliability assessment,
and operational market analysis \cite{BlHi09}. In general, it is
imperative to partition the grid judiciously in order to cope with
issues involving connected or disconnected ``subgrids.'' Regional
partitioning of the North American grid is illustrated in
Fig.~\ref{fig:nerc}, where each interconnection is further divided
into several zones for various planning and operation purposes.
However, the static and manual grid partitioning currently in
operation may soon become obsolete with the growing incorporation of
renewables and the overall system scaling.

The clustering criterion must be in accordance with grid
partitioning goals. In islanding applications, sub-groups of
generators are traditionally formed by minimizing the real
generator-load imbalance to regulate the system frequency within
each island. Recently, reactive power balance has been incorporated
in a multi-objective grid partitioning problem to support voltage
stability in islanding \cite{LiLi10}. For these methods, it is
necessary to reflect the real-time operating conditions that depend
on the slow-coherency among generators, and the flow density along
transmission lines.

Different from the islanding methods that deal with real-time
contingencies, zonal analysis intends to address the long-term
planning of transmission systems. Therefore, it is critical to
define appropriate distance metrics between buses. Most existing
works on long-term reliability have focused on the knowledge of
network topology, including the seminal work of \cite{WaSt98}, which
pointed out the ``small-world'' effects in power networks. To
account for the structure imposed by Kirchhoff's laws, it was
proposed in \cite{BlHi09} to define ``electrical distances'' between
buses using the inverse admittance matrix.

\section{Optimal Grid Operation}\label{sec:control}

Leveraging the extensive monitoring and learning modalities outlined in the previous section, the next-generation grid will be operated with significantly improved efficiency and reduced margins. After reviewing classical results on optimal grid dispatch, this section outlines challenges and opportunities related to demand-response programs, electric vehicle charging, and the integration of renewable energy sources with particular emphasis on the common optimization tools engaged.

\subsection{Economic Operation of Power Systems}\label{subsec:opf}

\subsubsection{Economic Dispatch}
Economic dispatch (ED) amounts to optimally setting the generation
output in an electric power network so that the load is served and
the cost of generation is minimized. ED pertains to generators which
consume some sort of non-renewable fuel in order to produce electric
energy, the most typical fuel types being oil, coal, natural gas, or
uranium. In what follows, a prototype ED problem is described, with
focus placed on a specific time span, e.g. 10 minutes or one hour,
over which the generation output is supposed to be roughly constant.

Specifically, consider a network with $N_g$ generators. Let $P_{G_i}$ be the output of the $i$th generator in MWh. The cost of the $i$th generator is determined by a function $C_i(P_{G_i})$, which represents the cost in \$ for producing energy of $P_{G_i}$ MWh (i.e., maintaining power output $P_{G_i}$ MW for one hour). The cost $C_i(P_{G_i})$ is modeled as strictly increasing and convex, with typical choices including piecewise linear or smooth quadratic functions. The output of each generator is an optimization variable
in ED, constrained within minimum and maximum bounds, $P_{G_i}^{\min}$ and $P_{G_i}^{\max}$, determined by the generator's physical characteristics~\cite[Ch.~2]{WollenbergBook}. Since once a power plant is on, it has substantial power output, $P_{G_i}^{\min}$ is commonly around 25\% of $P_{G_i}^{\max}$.

With $P_L$ denoting the load forecasted as described in Section \ref{subsubsec:load_prediction}, the prototype ED problem
is to minimize the total generation cost so that there is supply-demand
balance within the generators' physical limits:
\begin{subequations}
\label{ED}
\begin{align}
\min_{\{P_{G_i}\}}~~  & \sum_{i=1}^{N_g}  C_i(P_{G_i})
\label{ED-cost} \\
\text{subj.~to}~~ &  \sum_{i=1}^{N_g}  P_{G_i} = P_L \label{ED-balance}\\
                            &  P_{G_i}^{\min} \leq P_{G_i} \leq P_{G_i}^{\max}. \label{ED-lim}
\end{align}
\end{subequations}

Problem \eqref{ED} is convex, so long as the functions $C_i(P_{G_i})$ are convex. In this case, it can be solved very
efficiently. Convex choices of $C_i(P_{G_i})$ offer a model
approximating the true generation cost quite well and are used
widely in the literature. Nevertheless, the true cost in practice
may not be strictly increasing or convex, while the power output may
be constrained to lie in a collection of disjoint subintervals
$[P_{G_i}^{\min}, P_{G_i}^{\max}]$. These specifications make ED
nonconvex, and hence hard to solve. A gamut of approaches for
solving the ED problem can be found in \cite[Ch.~3]{WollenbergBook}.

Following a duality approach, suppose that Lagrange multiplier
$\lambda$ corresponds to constraint \eqref{ED-balance}. The
multiplier has units \$/MWh, which has the meaning of price. Then,
the KKT optimality condition implies that for the optimal generation
output $P_{G_i}^*$ and the optimal multiplier $\lambda^*$, it holds
that
\begin{equation}\label{netcosti}
P_{G_i}^* =\argmin_{P_{G_i}^{\min} \leq P_{G_i} \leq P_{G_i}^{\max}}
\{C_i(P_{G_i}) -\lambda^* P_{G_i} \}, \: i=1,\ldots,N.
\end{equation}
Due to \eqref{netcosti}, $C_i(P_{G_i}^*)$ is the $i$-th generator's cost in dollars. Moreover, if $\lambda^*$ is the price at which each generator is getting paid to produce
electricity, then $\lambda^*P_{G_i}^*$ is the profit for the $i$-th generator. Hence, the minimum in \eqref{netcosti} is the net cost, i.e., the cost minus the profit, for generator $i$. The latter reveals that the optimal generation dispatch is the one minimizing the net cost for each generator. If an electricity market is in place, ED is solved by the ISO, with $\{C_i(P_{G_i})\}$ representing the supply bids.

There are two take-home messages here. First, a very important operational feature of an electrical power network is to balance supply and demand in the most economical manner, and this can be cast as an optimization problem. Second, the Lagrange multiplier corresponding to the supply-demand balance equation can be readily interpreted as a price. However, the formulation in \eqref{ED} entails two simplifying assumptions: (i) it does not account for the transmission network; and (ii) it only pertains to a specific time interval, e.g., one hour. In practice, the power output across consecutive time intervals is limited by the generator physical characteristics. Even though the more complex formulations presented next alleviate these simplifications, the two take-home messages are still largely valid.

\subsubsection{Optimal Power Flow}
The first generalization is to include the transmission network, using the DC load flow model of Sec.~\ref{sec:modeling}; cf. \eqref{eq:pBtheta}. The resultant formulation constitutes the \emph{DC optimal power flow} (DC OPF) problem~\cite{ChWoWa00}. Specifically, it is postulated that at each bus there exist a generator and a load with output $P_{G_m}$, and demand $P_{L_m}$, respectively. The cases of no or multiple generators/loads on a bus can be readily accommodated.

Recall from~\eqref{eq:power_flow_dc_active} that the real power flow from bus $m$ to $n$ is approximated by $P_{mn} \approx -b_{mn} (\theta_m-\theta_n)$. The bus angles $\{\theta_m\}$ are also variables in the DC OPF problem that reads
\begin{subequations}
\label{DCOPF}
\begin{align}
&\min_{\{P_{G_m},\theta_m\}}~~   \sum_{m=1}^{N_b}  C_m(P_{G_m})
\label{DCOPF-cost} \\
&\text{subj.~to}~~  \notag \\
& P_{G_m} - P_{L_m}  = - \sum_{n\in\ccalN_m} 
b_{mn} (\theta_m-\theta_n), \: m=1,\ldots,N_b
\label{DCOPF-balance}\\
&   P_{G_m}^{\min} \leq P_{G_m} \leq P_{G_m}^{\max}, \: m=1,\ldots,N_b
\label{DCOPF-genlim} \\
&|P_{mn}| =|b_{mn} (\theta_m-\theta_n)| \leq P_{mn}^{\max}, \: m,n=1,\ldots,N_b.
\label{DCOPF-Txlim}
\end{align}
\end{subequations}
The objective in \eqref{DCOPF-cost} is the total generation cost.
Constraint \eqref{DCOPF-balance} is the per bus balance.
Specifically, the left-hand side of~\eqref{DCOPF-balance} amounts to
the net power injected to bus $m$ from the generator and the load
situated at the bus, while the right-hand side is the total power
that flows towards all neighboring buses. Upon defining vectors for
the generator and the load powers,~\eqref{DCOPF-balance} could be
written in vector form as
$\mathbf{p}_G-\mathbf{p}_L=\mathbf{B}_x\boldsymbol{\theta}$
[cf.~\eqref{eq:pBtheta}]. Finally, constraint \eqref{DCOPF-Txlim}
enforces power flow limits for line protection.

For convex generation costs $C_m(P_{G_m})$, the DC-OPF problem is convex too, and hence, efficiently solvable. A major consequence of considering per bus balance equations is that every bus may have a different Lagrange multiplier. The pricing interpretation of Lagrange multipliers implies that a different price, called locational marginal price, corresponds to each bus. The ED problem \eqref{ED} can be thought of as a special case of DC OPF, where the entire network consists of a single bus on which all generators and loads reside.

Due to the DC load flow approximation, the accuracy of the DC OPF greatly depends on how well assumptions (A1)-(A3) hold for the actual power system. For better consistency with (A2), it is further suggested to penalize the cost \eqref{DCOPF-cost} with the sum of squared voltage angle differences $\sum_{\mathrm{lines}} (\theta_m-\theta_n)^2$, which retains convexity. Even if the DC OPF is a rather simplified model for actual power systems, it is worth stressing that it is used for the day-to-day operation in several North American ISOs.

Consider next replacing the DC with the AC load flow model (cf. Sec.~\ref{sec:modeling}) in the OPF context. Generators and loads are now characterized not only by their real powers, but also the reactive ones, denoted as $Q_{G_m}$ and $Q_{L_m}$. The AC OPF takes the form
\begin{subequations}
\label{OPF}
\begin{align}
&\min_{\{P_{G_m},Q_{G_m},\mathcal{V}_m\}}~~   \sum_{m=1}^{N_b}  C_m(P_{G_m})
\label{OPF-cost} \\
&\text{subj.~to}~~ \notag \\
&P_{G_m} - P_{L_m} =  \sum_{n\in\ccalN_m} \mathrm{Re}\{\mathcal{S}_{mn}\} \notag \\
&Q_{G_m} - Q_{L_m} =  \sum_{n\in\ccalN_m} \mathrm{Im}\{\mathcal{S}_{mn}\}
\label{OPF-balance} \\
& \eqref{eq:smn},\quad \eqref{eq:i_mn}
\notag\\
&                 P_{G_m}^{\min} \leq P_{G_m} \leq P_{G_m}^{\max}; \; Q_{G_m}^{\min} \leq Q_{G_m} \leq Q_{G_m}^{\max} \label{OPF-genlim}\\
&|\mathrm{Re}\{\mathcal{S}_{mn}\}| \leq P_{mn}^{\max}; \;|\ccalS_{mn} | \leq S_{mn}^{\max}; \;
V_m^{\min}\leq |\mathcal{V}_m| \leq V_m^{\max}. \label{OPF-Txlim}
\end{align}
\end{subequations}
Constraint \eqref{OPF-balance} reveals that now both the real and
reactive powers must be balanced per bus. Recall further that
$\mathcal{S}_{mn}$ represents the complex power flowing over line
$(m,n)$. Therefore, the first constraint in \eqref{OPF-Txlim} refers
to the real power flowing over line $(m,n)$
[cf.~\eqref{DCOPF-Txlim}], while the second to the apparent power.
The last constraint in \eqref{OPF-Txlim} calls for voltage amplitude
limits.

Due to the nonlinear (quadratic equality) couplings between the
power quantities and the complex voltage phasors, the AC OPF in
\eqref{OPF} is highly nonconvex. Various nonlinear programming
algorithms have been applied for solving it, including the gradient
method, Newton-Raphson, linear programming, and interior-point
algorithms; see e.g., \cite[Ch. 13]{WollenbergBook}. These
algorithms are based on the KKT necessary conditions for optimality,
and can only guarantee convergence to a stationary point at best.
Taking advantage of the quadratic relations from voltage phasors to
all power quantities as in SE, the SDR technique has been
successfully applied, while a zero duality gap has been observed for
many practical instances of the AC OPF, and theoretically
established for tree networks; see \cite{Lavaei12,LZD12}, and
references therein. SDR-based solvers for three-phase OPF in distribution networks is considered in \cite{DaGiWo12}.

The AC OPF offers the most detailed and accurate model of the
transmission network. Two main advantages over its DC counterpart
are: i) the ability to capture ohmic losses; and ii) its flexibility
to incorporate voltage constraints. The former is possible because
the resistive part of the line $\pi$-model is included in the
formulation. Recall in contrast that assumption (A1) in the DC model
sets $r_{mn}=0$. But it is exactly the resistive nature of the line
that causes the losses. In view of \eqref{OPF}, the total ohmic
losses can be expressed as $\sum_{m} (P_{G_m}-P_{L_m})$. Such losses
in the transmission network may be as high as 5\% of the total load
so that they cannot be neglected \cite[Sec.~5.2]{ExpConCanBook}.

The discussion on OPF---with DC or AC power flow---so far has focused on economic operation objectives. System reliability is another important consideration, and the OPF can be modified in order to incorporate security constraints too, leading to the \emph{security-constrained OPF} (SCOPF). Security constraints aim to ensure that if a system component fails---e.g., if a line outage occurs---then the remaining system remains operational. Such failures are called contingencies. Specifically, the SCOPF aims to find an operating point such that even if a line outage occurs, all post-contingency system variables (powers, line flows, bus voltages, etc.) are within limits. The primary concern is to avoid cascading failures that are the main reasons for system blackouts. As explained in Sec.~\ref{subsec:topology}, if a line is in outage, the power flows on all other lines are adjusted automatically to carry the generated power.

SCOPF is a challenging problem due to the large number of possible contingencies. For the case of the DC OPF, power flows after a line outage are linearly related to the flows before the outage through the line outage distribution factors (LODFs) \cite{ChWoWa00}, \cite[Ch.~11]{WollenbergBook}. The LODFs can be efficiently calculated based on the bus admittance matrix $\mathbf{B}_x$ and are instrumental in the security-constrained DC OPF. The case of AC OPF is much more challenging, and a possible
approach is enumeration of all possible contingency cases; see e.g., \cite[Sec.~13.5]{WollenbergBook} for different approaches.

\subsubsection{Unit Commitment}
Here, the scope of DC OPF is broadened to incorporate the scheduling
of generators across multiple time periods, leading to the so-termed
\emph{unit commitment} (UC) problem. It is postulated that the
scheduling horizon consists of periods labeled as $1\ldots,T$ (e.g.,
a day consisting of 24 1-hour periods). Let $P_{G_m}^t$ be the
output of the $m$-th generator at period $t$, and $P_{L_m}^t$ the
respective demand. The generation cost is allowed to be
time-varying, and is denoted by $C_m^t(P_{G_m}^t)$. A binary
variable $u_m^t$ per generator and period is introduced, so that
$u_m^t=1$ if generator $m$ is on at $t$, and $u_m^t=0$ otherwise.
Moreover, the $m$th bus angle at $t$ is denoted by $\theta_m^t$.

Consideration of multiple time periods allows inclusion of practical
generator constraints into the scheduling problem. These are the
ramp-up/down and minimum up/down time constraints. The former
indicate that the difference in power generation between two
successive periods is bounded. The latter mean that if a unit is
turned on, it must stay on for a minimum number of hours; similarly,
if it is turned off, it cannot be turned back on before a number of
periods. The UC problem is formulated as follows.
\begin{subequations}
\label{UC}
\begin{align}
&\min_{\{P^t_{G_m},\theta_m^t,u_m^t\}}~~   \sum_{t=1}^T
\sum_{m=1}^{N_b}  \left[ C_m^t(P_{G_m}^t) +
S_m^t(\{u_m^\tau\}_{\tau=0}^t) \right]
\label{UC-cost} \\
&\text{subj.~to}~~~~ \notag \\
&  P_{G_m}^t = P_{L_m}^t - \sum_{n\in\ccalN_m}
b_{mn} (\theta_m^t-\theta_n^t), \notag\\
& \mspace{170mu} m=1,\ldots,N_b,\, t=1,\ldots,T
\label{UC-balance}\\
&  u_m^t P_{G_m}^{\min} \leq P_{G_m}^t \leq u_m^t P_{G_m}^{\max}, \; m=1,\ldots,N_b,\, t=1,\ldots,T
\label{UC-genlim} \\
& P_{G_m}^{t}-P_{G_m}^{t-1}\leq R_m^{\mathrm{up}};\, P_{G_m}^{t-1}-P_{G_m}^{t}
\leq R_m^{\mathrm{down}}, \notag \\
& \mspace{170mu} m=1,\ldots,N_b, \, t=1,\ldots,T
\label{UC-ramp}\\
& u_m^t - u_m^{t-1} \leq u_m^\tau, \; \tau=t+1,\ldots,\min\{t+T_m^{\mathrm{up}}-1,T\} ,\notag\\
& \mspace{290mu} t=2,\ldots,T 
\label{UC-uptime}\\
& u_m^{t-1} - u_m^t \leq 1- u_m^\tau, \notag \\
& \tau=t+1,\ldots,\min\{t+T_m^{\mathrm{down}}-1,T\} ,\, t=2,\ldots,T
\label{UC-downtime}\\
 &  |b_{mn} (\theta_m^t-\theta_n^t)| \leq P_{mn}^{\max} , \; m,n=1,\ldots,N_b,\, t=1,\ldots,T
 \label{UC-Txlim}\\
& u_m^t\in\{0,1\}, \; m=1,\ldots,N_b,\, t=1,\ldots,T.
\label{UC-bin}
\end{align}
\end{subequations}
The term $S_m^t(\{u_m^\tau\}_{\tau=0}^t)$ in the
cost~\eqref{UC-cost} captures generator start-up or shut-down costs.
Such costs are generally dependent on the previous on/off activity.
For instance, the more time a generator has been off, the more
expensive it may be to bring it on again. The initial condition
$u_m^0$ is known. It is also assumed that $C_m^t(0)=0$. The balance
equation is given next by \eqref{UC-balance}. Generation limits are
captured by \eqref{UC-genlim}. Constraint \eqref{UC-ramp} represents
the ramp-up/down limits, where the bounds $R_m^{\mathrm{up}}$ and
$R_m^{\mathrm{down}}$ and the initial condition $P_{G_m}^0$ are
given. Constraint \eqref{UC-uptime} means that if generator $m$ is
turned on at period $t$, it must remain on for the next
$T_m^{\mathrm{up}}$ periods; and similarly for the minimum down time
constraint in~\eqref{UC-downtime}, where both $T_m^{\mathrm{up}}$
and $T_m^{\mathrm{down}}$ are given~\cite{TaBi00}. The line flow
constraints are given by~\eqref{UC-Txlim}, while the binary feasible
set for the scheduling variables $u_m^t$ is shown in~\eqref{UC-bin}.

\begin{figure}
\centering
\includegraphics[width=0.80\linewidth]{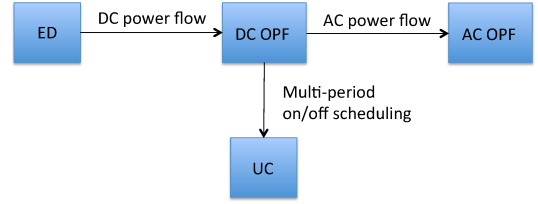}
\vspace*{-1em}
\caption{Relationship between the ED, OPF (DC and AC) and UC. From left to right, increasing detail in the transmission network model. From top to bottom, single- to multi-period scheduling (also applicable to ED and AC OPF).}
\label{fig:econ_opt}
\end{figure}

It is clear that problem~\eqref{UC} is a mixed integer program. What
makes it particularly hard to solve is the coupling across the
binary variables expressed by~\eqref{UC-uptime}
and~\eqref{UC-downtime}. Note that the DC OPF in~\eqref{DCOPF} is a
special case of the UC~\eqref{UC} with the on/off scheduling fixed
and the time horizon limited to a single period. It is noted in
passing that a multi-period version of the DC OPF can also be
considered, by adding the ramp constraints to~\eqref{DCOPF} while
keeping the on/off scheduling fixed in~\eqref{UC}, therefore
obtaining a convex program. Most importantly, note that the UC
dimension can be brought into the remaining two problems described
here, that is, the ED and the AC OPF. In the latter, the problem has
two mathematical reasons for being hard, namely, the integer
variables and the nonconvexity due to the AC load flow. The problems
discussed here are illustrated in Fig.~\ref{fig:econ_opt}.

A traditional approach to solving the UC is to apply Langrangian
relaxation with respect to the balance equations
\cite[Ch.~5]{WollenbergBook}, \cite{BeLaSaPo83}, \cite{TaBi00}. The
dual problem can be solved by a non-differentiable optimization
method (e.g., a subgradient or bundle method), while the Lagrangian
minimization step is solved via dynamic programming. An interesting
result within the Lagrangian duality framework is that the duality
gap of the UC problem without a transmission network diminishes as
the number of generators increases \cite{BeLaSaPo83}. One of the
state-of-the-art methods for UC is Benders decomposition, which
decomposes the problem into a master problem and tractable
subproblems \cite[Ch.~8]{ShYaLi02}.

\subsection{Demand Response}\label{subsec:dr}
Demand response (DR) or load response is the adaptation of end-user power
consumption to time-varying (or time-based) energy pricing, which is
judiciously controlled by the utility companies to elicit desirable
energy usage \cite{wollen11acc}, \cite{Hamilton-dr-pmag}. The smart grid
vision entails engaging residential end-users in DR programs. Residential
loads have the potential to offer considerable gains in terms of flexible
load response, because their consumption can be adjusted---e.g., an air
conditioning unit (A/C)---or deferred for later or shifted to an earlier
time. Examples of flexible loads include pool pumps or plug-in (hybrid)
electric vehicles. The advent of smart grid technologies have also made
available at the residential level energy storage devices (batteries), which
can be charged and discharged according to residential needs, and thus
constitute an additional device for control.

Widespread adoption of DR programs can bring significant benefits to the
future grid. First, the peak demand is reduced as a result of the load
shifting capability, which can have major economical benefits. Without DR,
the peak demand must be satisfied by generation units such as gas turbines
that can turn on and be brought in very fast during those peaks. Such units
are very costly to operate, and markedly increase the electricity wholesale
prices. This can be explained in a simple manner by recalling the ED problem
and specifically~\eqref{netcosti}. Considering a gas turbine that is brought
in and does not operate at its limits,~\eqref{netcosti} implies that
$\lambda^*=C'(P_{G_\mathrm{turbine}}^*)$. Expensive units have exactly very
high derivative $C'$, that is, increasing their power output requires a lot
of fuel.

\begin{figure}
\centering
\includegraphics[scale=0.43]{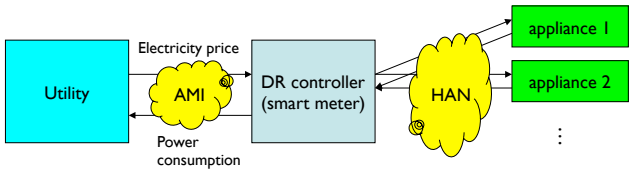}
\caption{Communications infrastructure facilitating DR capabilities.} \label{fig:ami-han}
\end{figure}

A second benefit of DR is that it has the potential to reduce the
end-user bills. This is due to the time-based pricing schemes, which
encourage consumption during reduced-price hours, but also because
the wholesale prices become less volatile as explained earlier,
which means that the electricity retailers can procure cheaper
sources. A third benefit is that DR can strengthen the adoption of
renewable energy. The reason is that the random and intermittent
nature of renewable energy can be compensated by the ability of the
load to follow such effects. More light into the latter concept will
be shed in Sec.~\ref{subsec:renewables}.

DR is facilitated by deployment of the \emph{advanced metering
infrastructure} (AMI), which comprises a two-way communication
network between utility companies and end-users (see
Fig.~\ref{fig:ami-han}) \cite{wollen11acc,Hamilton-dr-pmag}. Smart
meters installed at end-users' premises are the AMI terminals at the
end-users' side. These measure not just the total power consumption,
but also the power consumption profile throughout the day, and
report it to the utility company at regular time intervals---e.g.,
every ten minutes or every hour. The utility company sends pricing
signals to the smart meters through the AMI, for the smart meters to
adjust the power consumption profile of the various residential
electric devices, in order to minimize the electricity bill and
maximize the end-user satisfaction. Energy consumption is thus
scheduled through the smart meter. The communication network at the
customer's premises between the smart meter and the smart
appliances' controllers is part of the so-called home area network
(HAN).

Time-varying pricing has been a classical research
topic~\cite{CaBoSc82}. The innovation DR brings is that the
end-users' power consumption becomes controllable, and therefore,
\emph{part of the system optimization}. Novel formulations
addressing the various research issues are therefore called for.
DR-related research issues can be classified in two groups. The
first group deals with joint optimization of DR for a set of
end-users, which will be termed hereafter multi-user DR. The second
group focuses on optimal algorithm design for a single smart meter
with the aim of minimizing the electricity bill and the user
discomfort in response to real-time pricing signals. Each approach
has unique characteristics, as explained next.

Multi-user DR sets a system-wide performance objective accounting for the
cost of the energy provider and the user satisfaction. Joint scheduling must
be performed in a distributed fashion, and much of the effort is to come up
with pricing schemes that achieve this goal. Privacy of the customers must be
protected, in the sense that they do not reveal their individual power
consumption preferences to the utility, but the desired power consumption
profile is elicited by the pricing signals. One of the chief advantages of
joint DR scheduling for multiple users is that the peak power consumption is
reduced as compared to a baseline non-DR approach. The reason is that joint
scheduling opens up the possibility of loads being arranged across time so
that valleys are filled and peaks are shaved.

On the other hand, energy consumption scheduling formulations for a
single user can model in great detail the various smart appliance
characteristics, often leading to difficult nonconvex optimization
problems. This is in contrast with the vast majority of multi-user
algorithms, which tend to adopt a more abstract and less refined
description of the end-users' scheduling capabilities. More details
on the two groups of problems are given next.

\subsubsection{Multi-user DR}
\begin{figure}
\centering
\includegraphics[scale=0.5]{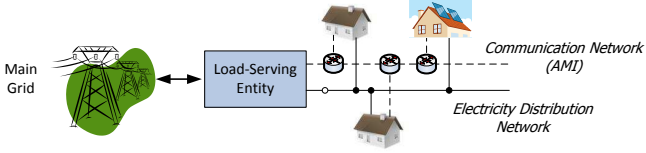}
\vspace*{-1em}
\caption{Power network consisting of electricity end-users and the
LSE.} \label{fig:lse}
\end{figure}

Consider $R$ residential end-users, connected to a single
load-serving entity (LSE), as illustrated in Fig~\ref{fig:lse}.  The
LSE can be an electricity retailer or an aggregator, whose role is
to coordinate the $R$ users' consumption and present it as a larger
flexible load to the main grid. The time horizon consists of $T$
periods, which can be a bunch of 1-hour or ten-minute intervals.
User $r$ has a set of smart appliances $\mathcal{A}_r$. Let
$p_{ra}^t$ be the power consumption of appliance $a$ of user $r$ at
time period $t$ (typically in kWh), and $\mathbf{p}_{ra}$ a $T
\times 1$ vector collecting the corresponding power consumptions
across slots.

The LSE incurs cost $C^t(s^t)$ for providing energy $s^t$ to the
users. This cost is essentially the cost of energy procurement from
the wholesale market or through direct contracts with energy
generation units,  and may also include other operation and
maintenance costs. Each user also adopts a utility function
$U_{ra}(\mathbf{p}_{ra})$, which represents user willingness to
consume power.

The prototype multi-user DR problem takes the following form
\begin{subequations}
\label{MU-DR}
\begin{align}
\min_{\{s^t\},\{p_{ra}^t\}} ~~& \sum_{t=1}^T C^t(s^t) - \sum_{r=1}^R
\sum_{a\in\mathcal{A}_r}  U_{ra}(\mathbf{p}_{ra})
\label{MU-DR-obj}\\
\text{subj.~to}~ &  s^t=\sum_{r=1}^R \sum_{a\in\mathcal{A}_r} p_{ar}^t,
\; ~~t=1,\ldots,T
\label{MU-DR-balance}\\
&  \mathbf{p}_{ra} \in \mathcal{P}_{ra}, \; a\in\mathcal{A}_r, ~~r=1,\ldots,R
\label{MU-DR-applim}\\
& s^{\min} \leq s^t \leq s^{\max}, \;~~ t=1,\ldots,T.
\label{MU-DR-lselim}
\end{align}
\end{subequations}
Clearly, the objective is optimizing the system's \emph{social welfare}.
Constraint~\eqref{MU-DR-balance} amounts to a balance equation for each
period. Moreover, the set $\mathcal{P}_{ra}$ in~\eqref{MU-DR-applim}
represents the scheduling constraints for every appliance, while
constraint~\eqref{MU-DR-lselim} bounds the power provided by the LSE.

Problem~\eqref{MU-DR} is convex as long as $C^t(s^t)$ is convex,
$U_{ra}(\mathbf{p}_{ra})$ is concave, and sets $\mathcal{P}_{ra}$
are convex. This is typically the case, and different works in the
literature address DR using versions of the previous formulation
\cite{ChLiJiLo12, Rad-TSG10, SaScWo11, GaGG-TSG12}. Various examples
of appliance models---including batteries---together with their
utility functions and constraint sets can also be found in the
aforementioned works.

Problem~\eqref{MU-DR} as described so far amounts to \emph{energy consumption
scheduling}. Another instance of DR that can be described by the previous
formulation is \emph{load curtailment}. In this context, there is an energy
deficit in the main grid for a particular time period, and the LSE must
regulate the power consumption to cover for this deficit. The situation can
be captured in~\eqref{MU-DR} by setting $T=1$ (single time period), and the
power deficit as $s^{\min} = s^{\max} = s$.  The cost $C^t$ does not affect
the optimization, while the negative of $U_{ra}(p_{ra})$ represents the
discomfort of the end-user due to the power curtailment, so the total
discomfort $-\sum_{r,a}U_{ra}(p_{ra})$ is minimized. This problem is
addressed in~\cite{PaHiGr10, KiKiYaTh11} and the references therein.

One of the main research objectives regarding \eqref{MU-DR} is to solve the scheduling problem in a distributed fashion, without having the functions $U_{ra}(\mathbf{p}_{ra})$ and sets $\mathcal{P}_{ra}$ communicated to the LSE in order to respect customer privacy. Algorithmic approaches typically entail message exchanges between the LSE and the users or among the users, and lead to different pricing interpretations and models. Specific approaches include gradient projection \cite{ChLiJiLo12}; block coordinate descent \cite{Rad-TSG10}; dual decomposition and subgradient method \cite{GaGG-TSG12}, \cite{PaHiGr10}; the Vickrey-Clarke-Groves (VCG) mechanism  \cite{SaScWo11}; Lagrangian relaxation and Newton method \cite{PaHiGr10}; and dual decomposition with the bisection and Illinois methods \cite{KiKiYaTh11}.

%

Formulation~\eqref{MU-DR} refers to ahead-of-time scheduling. Real-time
scheduling is also important. A real-time load response approach
operating on a second-to-second scale is developed in~\cite{KoMaCa11} and
references thereof. The aim is to have the aggregate power consumption of a
set of thermostatically controlled loads (TCLs), such as A/C units, follow a
desired signal. Model predictive control is employed to this end. Moreover,
in order to come up with a simple description of the state space model
pertaining to the set of TCLs, system identification ideas are brought to
bear.

\subsubsection{Single-user DR}
The problems here focus on minimizing the total cost due to energy
consumption or the peak instantaneous cost over a billing interval
(or possibly a combination thereof). User comfort levels and
preferences must also be taken into account.

Detailed modeling of appliance characteristics and scheduling
capabilities typically introduces integer variables into the
formulation, which is somewhat reminiscent of the unit commitment
problem [cf.~\eqref{UC}]; see e.g.,~\cite{Pedrasa-ResDER,
SJGG-CAMSAP11, SoWeSaJo11} and references therein. Solution
approaches include standard mixed-integer programming
techniques---e.g., branch-and-bound, Lagrangian relaxation, dynamic
programming---as well as random search methods such as genetic
algorithms and particle swarm optimization. An interesting result is
that when the problem is formulated over a continuous time horizon
and accounts for the fact that appliances can be turned on or off
anytime within the horizon, then it has zero duality
gap~\cite{GaGG-CDC11}.

Real-time approaches have also been pursued. A linear programming DR model
with robustness against price uncertainty and time-series-based
price prediction from period to period is developed
in~\cite{Conejo-RTDR}. Moreover, \cite{NgBeMaPa11} focuses on TCLs,
and specifically, on a building with multiple zones, with each zone
having its own heater. The aim is to minimize the peak instantaneous
cost due to the power consumption of all heaters, while keeping each
zone at a specified temperature interval. The problem is tackled
through a decomposition into a master mixed-integer program and per
zone heater control subproblems.

\subsection{Plug-in (Hybrid) Electric Vehicles}\label{subsec:phev}
As an important component of the future smart grid vision, electric
vehicles (EVs) including plug-in (hybrid) EVs (P(H)EVs) are
receiving a lot of attention. A global driving factor behind the
research and development efforts on EVs is the environmental concern
of the greenhouse gases emitted by the conventional fossil
fuel-based transportation. As the future grids accommodate the
renewable energy resources in an increasing scale, the carbon
footprint is expected to be markedly curbed by high EV penetration.
Electric driving also bears strategic relevance in the context of
growing international tension over key natural resources including
crude oil. From the simple perspective of improving overall energy
efficiency, electrification of transportation offers an excellent
potential.

PEVs interact directly with the power grid through plug-in charging of
built-in batteries. As such, judicious control and optimization of PEV
charging pose paramount challenges and opportunities for the grid economy and
efficiency. Since PEV charging constitutes an elastic energy load that can be
time-shifted and warped, the benefits of DR are to be magnified when PEV
charging is included in DR programs. In fact, as the scale of PEV adoption
grows, it is clear that smart coordination of the charging task will become
crucial to mitigate overloading of current distribution
networks~\cite{ClN10,WuA11,DMM11}. Without proper coordination, PEV charging
can potentially create new peaks in the load curves with detrimental effects
on generation cost. On the other hand, it is possible for the PEV aggregators
that have control over a fleet of PEVs to provide ancillary services by
modulating the charging rate in the vehicle-to-grid (V2G) concept
\cite{SoS12}. This in turn allows the utilities to depend less on
conventional generators with costly reserve capacities, and facilitates
mitigation of the volatility of renewable energy resources integrated to the
grid \cite{KeT05}. The aforementioned topics are discussed in more detail
next.

\subsubsection{Coordination of PEV Charging}
It is widely recognized that uncoordinated PEV charging can pose serious
issues on the economy of power generation and the quality of power delivered
through the distribution networks. PEVs are equipped with batteries with
sizable capacities, and it is not difficult to imagine that most people would
opt to start charging their vehicles immediately after their evening commute,
which is the time of the day that already exhibits a significant peak in
power demand~\cite{DMM11}. Fortunately, the smart grid AMI reviewed in Sec.~\ref{subsec:dr} provides the groundwork for effective scheduling and control of PEV charging to meet the challenges and sustain mass adoption.

A variety of approaches have been proposed for PEV charging coordination. The
power losses in the distribution network were minimized by optimizing the
day-ahead charging rate schedules for given PEV charging demands
in~\cite{ClN10}. Real-time coordination was considered in~\cite{DMM11}, where
the cost due to time-varying electricity price as well as the distribution
losses were minimized by performing a simple sensitivity analysis of the cost
and accommodating the charging priorities. Extending recent results on
globally optimal solution of the OPF problem via its Lagrangian dual
\cite{Lavaei12}, the optimality of similar approaches for PEV coordination
problems was investigated in \cite{SoL11}.

Interestingly, PEV charging can be also pursued in a distributed fashion.
Further, optimizing feeder losses of distribution networks, load factor, and
load variance are oftentimes equivalent problems~\cite{SHM11}. Leveraging the
latter, minimization of load variance was investigated in~\cite{GTL11}.
Specifically, the optimal day-ahead charging profiles ${\bf r}_n := [r_n(1),
\ldots, r_n(T)]$ for vehicle $n \in \{1,\ldots,N\}$ over a $T$-slot horizon,
are obtained by solving
\begin{subequations}\label{eq:valleyfill}
\begin{align}
\min_{{\bf r}_1,\ldots,{\bf r}_N} ~&~\sum_{t=1}^T \left(D(t)
+ \sum_{n=1}^N r_n(t) \right)^2 \label{eq:central_obj}\\
\textrm{subj. to }~&~ \underline{\bf r}_n \preceq {\bf r}_n
\preceq \overline{\bf r}_n, \quad n=1,\ldots,N\\
& \sum_{t=1}^T r_n(t) = B_n, \quad n=1,\ldots,N
\end{align}
\end{subequations}
where $D(t)$ is the given base demand, $\underline{\bf r}_n$ and
$\overline{\bf r}_n$ specify the limits on charging rates, and $B_n$
represents the total energy expended for charging PEV $n$ to the
desired state-of-charge~(SoC). The formulation is referred to as
``valley-filling" in~\cite{GTL11}, as it schedules PEV loads in the
valleys of the base load curve.

An optimal solution to~\eqref{eq:valleyfill} can be obtained iteratively
\cite{GTL11}. Supposing that the initial pricing signal $p^k(t) = D(t)$, $t =
1,2,\ldots,T$, and the initial charging profiles ${\bf r}_n^k(t)$ are
identically zero for iteration $k = 0$, each PEV $n$ updates charging
profiles ${\bf r}_n^{k+1}$ via
\begin{subequations}
\begin{align}
\min_{{\bf r}_n} ~&~\sum_{t=1}^T p^k(t) r_n(t) + \frac{N}{2} \left( r_n(t) - r_n^k(t) \right)^2 \label{eq:distr_obj}\\
\textrm{subj. to} ~&~ \underline{\bf r}_n \preceq {\bf r}_n \preceq \overline{\bf r}_n \textrm{ and} \sum_{t=1}^T r_n(t) = B_n.
\end{align}
\end{subequations}
A central entity such as the utility or a PEV aggregator then collects the profiles $\left\{{\bf r}_n^{k+1}\right\}$ from
all PEVs, and updates the pricing signal as
\begin{align}
p^{k+1}(t) = D(t) + \sum_{n=1}^N r_n^{k+1}(t). \label{eq:price_update}
\end{align}
The new pricing signals are then fed back to the PEVs and the procedure
iterates until convergence. It is clear from \eqref{eq:price_update} that the
per-vehicle objective in~\eqref{eq:distr_obj} corresponds to a first-order
estimate of the overall objective in~\eqref{eq:central_obj}, augmented with a
proximal term. The overall procedure turns out to be a projected gradient
search.

\subsubsection{Integration with Renewables and V2G}
It is only when the wide adoption of PEVs is coupled with large-scale
integration of renewable energy sources that the emission problem can be
alleviated, as the conventional generation itself contributes heavily to the
emission. However, renewable energy sources are by nature intermittent, and
often hard to predict accurately. By allowing the PEV batteries or fuel cells
to supply their stored power to the grid based on the V2G concept, it was
observed in~\cite{KeT05} that photovoltaic (PV) resources harnessed by the
EVs could competitively provide peak power (since the PV power becomes
highest few hours earlier than the daily load peak quite predictably), and
large-scale wind power could be stabilized for providing base power, via
intelligent control. For specific control strategies to accomplish such
benefits, formulations that maximize the profit for providing ancillary
services were considered in \cite{SoS12} and references therein.

\subsubsection{Charging Demand Prediction}
An important prerequisite task to support optimal coordination of
PEVs is modeling and prediction of the PEV charging demand. The
probability distributions of the charging demand were characterized
in \cite{LKP12} and references therein. Spatio-temporal PEV
charging demand was analyzed for highway traffic scenarios using a
fluid traffic model and a queuing model in~\cite{BaK12}. However,
there are many interesting issues remaining that deserve further
research in this forecasting task.

\subsection{Renewables}
\label{subsec:renewables}

The theme of Sec.~\ref{subsec:opf} has been economic scheduling of
generators, which consume non-renewable fuels. The subject of the
present section is on including generation from renewable energy
sources (RESs), with the two prime examples being wind and solar
energy. RESs are random and intermittent, which makes them
\emph{nondispatchable}. That is, RESs are not only hard to predict,
but their intermittency gives rise to high variability even within
time periods as short as 10 minutes. Therefore, they cannot be
readily treated as conventional generators, and be included in the
formulations of Sec.~\ref{subsec:opf}. In this context, methods for
integrating generation from RESs to the smart grid operations are
outlined next.

\subsubsection{Forecast-Based Methods}
\label{subsubsec:RES-forecast} To illustrate the forecast-based
methods, recall the ED problem [cf.~\eqref{ED}], and suppose that
there is also a wind power generator that can serve the load. The
output of the wind power generator for the next time period is a
random variable denoted by $W$. It is assumed that a forecast
$\hat{W}$ is available, and that the wind power generator has no
cost (as it does not consume fuel). Then, the balance constraint is
replaced by [cf.~\eqref{ED-balance}]
\begin{equation}
\label{eq:RES-forecast}
\sum_{i=1}^{N_g}  P_{G_i} = P_L - \hat{W}
\end{equation}
while the remainder of the ED problem remains the same. Since the
load is actually forecasted (cf.
Sec.~\ref{subsubsec:load_prediction}),
constraint~\eqref{eq:RES-forecast} essentially treats the uncertain
RES no different than a negative load.

In order for the forecast to be accurate, the time period of ED is
recommended to be short, such as 10 minutes. Building on this, a
multi-period ED is advocated in~\cite{Ilic-TPS11}, where the main
feature is a model-predictive control approach with a moving
horizon. Specifically, the ED over multiple periods and accompanying
forecasts is solved for e.g., 6 ten-minute periods representing an
hour. The generation is dispatched during the first period according
to the obtained solution. Then, the horizon is moved, and a new
multi-period ED with updated forecasts is solved, whose results are
applied only to the next period, and so on. Such a method can
accommodate the ramping constraints, and is computationally
efficient.

\subsubsection{Chance-Constrained Methods}
\label{susubsec:RES-chance} To account for the random nature of RES
in ED, the probability distribution of $W$ comes handy. Specifically, the
constraint is now that the supply-demand balance holds with high probability
$\varepsilon$, say 99\%. Hence, \eqref{eq:RES-forecast} is substituted by the chance constraint
\begin{equation}
\label{eq:RES-chance}
\mathrm{Prob}\left[\sum_{i=1}^{N_g}  P_{G_i} + W \geq P_L\right]\geq \varepsilon.
\end{equation}
Note that the equality of the balance equation has been replaced by an inequality in~\eqref{eq:RES-chance}, because excess power from RESs can in principle be curtailed.

To solve the chance-constrained ED, the distribution of $W$ must be known. For wind power, this is derived from the wind speed distribution, and the speed-power output mapping of the generator~\cite{Liu10}. The most typical speed distribution is Weibull, while the speed-power output mapping is nonlinear. Evidently, this approach poses formidable modeling and computing challenges when multiple RESs and their spatio-temporal correlation are considered. The probability that the load is not served (immediately obtained from the one in~\eqref{eq:RES-chance}) is often called loss of load probability. Related sophisticated methods which account for chance constraints are also described in~\cite{VaWuBi11}. An alternative approach not requiring the joint spatio-temporal wind distribution is presented in \cite{ZhGaGia13}.

\subsubsection{Robust (Minmax) Optimization}
\label{susubsec:RES-robust} This approach postulates that the power
generation from all RESs across space and time belongs to a
deterministic uncertainty set. The aim is to minimize the worst-case
operational costs, while setting the dispatchable generation and
other optimization variables to such levels so that the balance is
satisfied for any possible RES output within the uncertainty set.
The main attractive feature here is that no detailed probabilistic
models are needed. Only the uncertainty set must be obtained, e.g.,
from historical data, or, meteorological factors.

A robust version of UC [cf.~\eqref{UC}] is presented next. Following the
notation of Sec.~\ref{subsec:opf}, it is postulated that there are RESs with
power output $W^t_m$ per bus and time period. Let
$\mathbf{w}:=\{W_m^t\}_{m,t}$, and $\mathcal{W}$ denote the uncertainty set
for $\mathbf{w}$. The optimization variables are set in two stages. The
on/off variables $\mathbf{u}:=\{u_m^t\}_{m,t}$ are chosen during the first
stage. The power generation variables and bus angles are set after the RES
power output is realized---which constitutes the second stage. Therefore, the
power outputs and bus angles are functions of the commitments as well as the
RES power outputs, and are denoted as $P_{G_m}^t(\mathbf{u}, \mathbf{w})$ and
$\theta_m^t(\mathbf{u}, \mathbf{w})$. The robust two-stage UC problem takes
the form
%
\begin{subequations}
\label{RES-rob}
\begin{align}
&\min_{\mathbf{u}, \{P_{G_m}^t(\mathbf{u}, \mathbf{w}),\theta_m^t(\mathbf{u}, \mathbf{w})\} }~~   \sum_{t=1}^T \sum_{m=1}^{N_b}
S_m^t(\{u_m^\tau\}_{\tau=0}^t) \notag \\
& \mspace{150mu} + \max_{\mathbf{w}\in\mathcal{W}} \sum_{t=1}^T
\sum_{m=1}^{N_b}  C_m^t(P_{G_m}^t(\mathbf{u}, \mathbf{w}))
\label{RES-rob-obj}\\
& \text{subj.~to}~~ \notag\\
&  \eqref{UC-uptime}, \eqref{UC-downtime}, \eqref{UC-bin}
\label{RES-rob-1ststage}\\
& \hspace*{-0.7em}\left.\begin{array}{l}
\eqref{UC-genlim}, \eqref{UC-ramp}, \eqref{UC-Txlim}\\
P_{G_m}^t(\mathbf{u}, \mathbf{w}) + W^t_m = P_{L_m}^t \\
\mspace{150mu} + \sum\limits_{n\in\mathcal{N}_m}
b_{mn} [\theta_m^t(\mathbf{u}, \mathbf{w})-\theta_n^t(\mathbf{u}, \mathbf{w})]
\end{array}\right\} \notag\\
& \mspace{320mu}\forall \mathbf{w}\in\mathcal{W}.
\label{RES-rob-2ndstage}
\end{align}
\end{subequations}

The objective~\eqref{RES-rob-obj} consists of the startup/shutdown
costs related to the on/off scheduling decisions, as well as the
worst-case generation costs. The constraints
in~\eqref{RES-rob-1ststage} pertain only to the on/off variables,
and are identical to those in the UC problem. The remaining UC
constraints must be satisfied for all possible realizations of the
uncertain RES, as indicated in~\eqref{RES-rob-2ndstage}.

The solution of problem~\eqref{RES-rob} proceeds as follows. The
on/off decisions $u_m^t$ determine the UC ahead of the horizon
$\{1,\ldots,T\}$. Then, at each period, after the RES power output
is realized, functions $P_{G_m}^t(\mathbf{u}, \mathbf{w})$ yield the
power generation dispatch. The punch line of this two-stage robust
program is that generation becomes adaptive to the RES uncertainty.
Solution methods typically involve pertinent decompositions and
approximations~\cite{ZhaoZ10, BeLiSu10}.

A different robust approach for energy management in microgrids is pursued
in~\cite{YZNGGG-SGComm12}. Microgrids are power systems comprising many
distributed energy resources (DERs) and electricity end-users, all deployed
across a limited geographical area. Depending on their origin, DERs can come
either from distributed generation (DG), meaning small-scale power generators
based on fuels or RESs, or from distributed storage (DS), such as batteries.
The case where a microgrid is connected to the main grid, while energy can be
sold to or purchased from the main grid, is considered
in~\cite{YZNGGG-SGComm12}. The approach adopts a worst-case transaction cost.
Leveraging the dual decomposition, its solution is obtained in a distributed
fashion by local controllers of the DG units and dispatchable loads.

\begin{figure}
\centering
\includegraphics[width=1.0\linewidth]{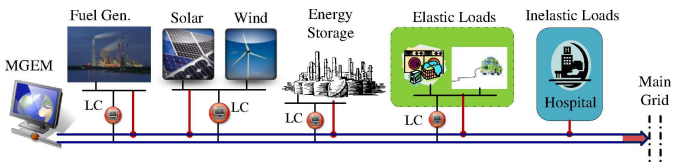}
\vspace*{-1em}
\caption{Distributed control and computation architecture of a
microgrid system. The microgrid energy manager (MGEM) coordinates
the local controllers (LCs) of DERs and dispatchable loads. }
\label{fig:MGModel}
\vspace*{-1em}
\end{figure}

\subsubsection{Scenario-based Stochastic Programming}
\label{susubsec:RES-scenario} This method also amounts to a
two-stage adaptive approach, albeit in a different manner than the
previous one. Here, a discrete set of possible scenarios for the RES
power output across the horizon is considered. For instance,
considering 8 hours with power output taking 7 possible values,
there are $7^8$ possible scenarios. A probability is attached to a
each of these scenarios (or only to a selection thereof). Similar
to~\eqref{RES-rob-obj}, the objective includes startup/shutdown
costs due to on/off scheduling. But instead of a worst-case part,
the expected cost of generation dispatch with respect to the
scenario probabilities is included in the objective.

The aforementioned approach is pursued in~\cite{BoGa08}, whereby the
scheduling of spinning reserves is also included. Spinning reserve is
generation capacity that is not currently used to serve the load, but is
connected to the system (spinning) and is available to serve the load in case
there is loss of generation. Spinning reserves are instrumental components of
any power system, and the premise here is that they can be provisioned in a
manner adaptive to the RES uncertainty.

\subsubsection{Multi-Stage Stochastic Dynamic Programming}
\label{susubsec:RES-DP} The aim here is to address the decision
making challenges for an LSE obtaining energy from the market as
well as from RESs (cf. Fig.~\ref{fig:lse}). The LSE may procure
energy in the day-ahead market, as well as in the real-time market,
which is a decision made on-the-fly during the scheduling horizon.
The energy from RESs is typically cost-free, but random. In
addition, the LSE must provide power to the end-users during the
horizon, and take the associated pricing decisions. The
multiple-timescale feature reflects exactly the resolution over
day-ahead, hour-ahead, and real-time (e.g., at the scale of minutes)
decisions. Multi-stage dynamic programming captures the coupling of
decisions across time due to end-users' power requirements---e.g.,
total energy requested over a specified interval, or, price-adaptive
random opportunistic demand \cite{JiLow-Allerton11, HeMuZh11}.

\subsubsection{Network Optimization Based on Long-Term Average
Criteria}\label{susubsec:RES-netopt} This approach relies on
queueing-theoretic and Lyapunov-based stochastic network optimization methods
popular in resource allocation tasks for wireless networks. A load-serving
entity obtaining energy from the market as well as from RESs is considered
in~\cite{Neely10, HuWaRa11}. The objective is cost minimization or social
welfare maximization in a long-term average fashion over an infinite horizon,
and the decision variables include pricing and power provided to end-users;
see also~\cite{KoHaTa11} for energy storage management policies.

\section{Open Issues}\label{sec:open}
Although the SP research efforts on power grid are fast growing, there are a lot of open issues awaiting investigation. Regarding situational awareness, integrating local power grids into interconnections poses modeling and computational challenges. Monitoring grids of dimensionality and detail calls for scalable and modular algorithms. To communicate and process the massive volume of measurements in real time with tractable complexity, the issues related to compressing, layering, relaying, and storing these data must be considered too. The ``big data'' challenges further extend to addressing the missing data and the under-determinacy of the resultant systems of equations, as well as model reduction tasks, for which contemporary statistical learning approaches could provide viable solutions.

The control and optimization dimensions entail conventional generation as well as RESs, interconnected via transmission and distribution networks, serving large industrial customers and residential end-users with smart appliances and P(H)EVs, as well as microgrids with distributed generation and storage. SP researchers can cross-fertilize their ample expertise on resource allocation gained in the context of communication networks to optimize power network operations. Major challenges include the successful coordination of system-level economic operations such as OPF and UC, while embracing small-scale end-users through DR and coordinated P(H)EV charging. Integrating random and intermittent RESs across all levels poses further challenges. Issues related to leveraging the markedly improved monitoring modalities in grid operations are worth careful study. Albeit research efforts tackling individual problems have yielded promising outcomes, achieving the grand goal of reliable and efficient grid operations still calls for novel formulations, insightful approximations, integration, and major algorithmic breakthroughs.

\bibliographystyle{IEEEtranS}
\bibliography{IEEEabrv,power_literature}

\end{document}